\documentclass[12pt]{article}  

\usepackage{amsmath,amsxtra,amssymb,latexsym, amscd}
\usepackage[mathscr]{eucal}

\begin{document}

\title{ \huge\bf Finite Hilbert Transforms\\
Logarithmic Potentials and Singular Integral Equations
 }

\def\1{\rule{0cm}{0cm}} \def\qd{\rule{3mm}{3mm}} \def\BB{$\bullet$}
\renewcommand{\arraystretch}{1.25}
\renewcommand{\theequation}{\thesectn.\arabic{equation}}
\def\sce{\setcounter{equation}{0}}  \newcounter{sectn} 
\newcounter{sbsect}
\def\sect#1{\addtocounter{section}{1}\sce\setcounter{sbsect}{0}%
        \renewcommand{\thesectn}{\thesection}\1\smallskip\\
        {\1\hspace{-2em}\large\bf\thesectn.\qquad #1\smallskip\par}}
\def\subsect#1{\addtocounter{sbsect}{1}\sce%
        
\renewcommand{\thesectn}{\thesection:\Alph{sbsect}}\1\smallskip\\
        {\bf\1\hspace{-1.5em}\thesectn.\qquad #1\smallskip\par}}
\newtheorem{Theorem}{THEOREM} \newtheorem{Lemma}[Theorem]{LEMMA}
\newtheorem{Corollary}[Theorem]{COROLLARY}
\def\thm#1#2{\be{Theorem}{\lb{#1} #2}} \def\LEM#1#2{\BE{Lemma}{\LB{#1} 
#2}}
\def\COR#1#2{\BE{Corollary}{\LB{#1} #2}}
\def\proof{\bigskip\noindent {\sc Proof:}\qquad}
\def\REM{\1\smallskip\par\noindent{\bf REMARK:}\qquad }
\def\qed{\hfill$\quad$\qd\medskip\\} \def\ds{\displaystyle}
\def\LB#1{\label{#1}} \def\BE#1#2{\begin{#1} #2 \end{#1}}
\def\EQ#1#2{\BE{equation}{\LB{#1} #2}} \def\ARR#1#2{\BE{array}{{#1} 
#2}}
\def\DES#1{\BE{description}{#1}} \def\QT#1{\BE{quote}{#1}}
\def\ENUM#1{\BE{enumerate}{#1}} \def\ITM#1{\BE{itemize}{#1}}
 \def\COM#1{\par\noindent{\bf COMMENT:\quad\sl #1}\par\noindent}
\def\mapsfrom{\hbox{$\;{\leftarrow}\kern-.15em{\mapstochar}\:\:$}}
\def\vv{\kern.344em{\rule[.18ex]{.075em}{1.32ex}}\kern-.344em}
\def\RE{\mbox{\rm I\kern-.21em R}} \def\CX{\mbox{\rm \vv C}}
\def\imp{\Rightarrow} \def\emb{\hookrightarrow} 
\def\wk{\rightharpoonup}
\def\rd{\dot{\1}} \def\d{\cdot} \def\+{\oplus} \def\x{\times}
\def\<{\langle} \def\>{\rangle} \def\o{\circ} \def\at#1{\Bigr|_{#1}}
\def\cd{\partial} \def\grad{\nabla} \def\L{\left} \def\R{\right}
\def\bx{\mathbf{x}} \def\by{\mathbf{y}} \def\bS{\mathbf{S}}
\def\I{{\cal I}} \def\A{{\cal A}} \def\D{{\cal D}}\def\bc{\mathbf{c}}
\def\bx{\mathbf{x}} \def\by{\mathbf{y}} \def\bS{\mathbf{S}}
\def\H{{\mathcal H}} \def\U{{\cal U}} \def\D{{\cal 
D}}\def\bc{\mathbf{c}}
\def\eq{equation} \def\de{differential \eq} \def\pde{partial \de}
\def\sol{solution} \def\pb{problem} \def\bdy{boundary} 
\def\fn{function}
\def\dde{delay \de} \def\ev{eigenvalue}
\def\R{\mathbb R}
\def\C{\mathbb C}
\author{
{\bf Dang Vu Giang}\\
Hanoi Institute of Mathematics\\
18 Hoang Quoc Viet, 10307 Hanoi, Vietnam\\
{\footnotesize          e-mail: $\<$dangvugiang@yahoo.com$\>$}\\
\1\\
}
\maketitle 
{\footnotesize
\noindent {\bf Abstract.} Several  interesting formulas concerning finite Hilbert transform and logarithmic integrals are proved with application in determining equilibrium measures, planar limits of analytic random matrix models with $1-$cut potential  and solving singular integral equations.

\medskip


\medskip

}

\bigskip

\sect{\bf  Hilbert Transforms and complex Hardy Spaces}

\bigskip
\par\noindent We define the Hilbert transform and the real Hardy space ${{H}^{1}}\left( \mathbb{R} \right).$ The Hilbert transform $Hf=\tilde{f}$ of a function $f\in {{L}^{p}}\left( \mathbb{R} \right)$    $\left( 1\le p<\infty  \right)$ is defined by letting
$$Hf\left( x \right)=\tilde{f}\left( x \right)=\frac{1}{\pi }\left( \text{p}\text{.v}\text{.} \right)\int\limits_{-\infty }^{\infty }{\frac{f\left( t \right)}{x-t}}\cdot dt.$$
For example, the Hilbert transform of the characteristic function  ${{\chi }_{\left( a,b \right)}}$ of the interval   $\left( a,b \right)$ is 
\[{{\tilde{\chi }}_{\left( a,b \right)}}\left( x \right)=\frac{1}{\pi }\cdot \ln \left| \frac{x-a}{x-b} \right|.\]
Moreover, for a tipical function 
\[\begin{aligned}
   f\left( x \right)&=\frac{1}{\sqrt{1-{{x}^{2}}}}\text{  for }\left| x \right|<1 \\ 
 & =0\quad\text{         otherwise}  
\end{aligned}\]
we have
\[\begin{aligned}
\tilde{f}\left( x \right)  & =\frac{1}{\pi }\int\limits_{-1}^{1}{\frac{1}{\sqrt{1-{{y}^{2}}}}\frac{dy}{x-y}}\text{ } \\ 
 & =\frac{2}{\pi }\int\limits_{0}^{\infty }{\frac{dt}{x\left( 1+{{t}^{2}} \right)-\left( 1-{{t}^{2}} \right)}}=0\quad\text{    if }\left| x \right|<1 \\ 
 & =\frac{1}{\sqrt{{{x}^{2}}-1}}\quad\text{    if }x>1 \\ 
 & =-\text{ }\frac{1}{\sqrt{{{x}^{2}}-1}}\quad\text{   if }x<-1.  
\end{aligned}\]
(Note that formula (8)  in page 175 of \cite{Tricomi} is incorrect.) 
Here, we change the variable
\[y=\frac{1-{{t}^{2}}}{1+{{t}^{2}}}\qquad \frac{dy}{\sqrt{1-{{y}^{2}}}}=-\frac {2dt}{1+{{t}^{2}}}\]
 and use the fact that for $\alpha,\beta >0$
\[\begin{aligned}
  & \int\limits_{0}^{\infty }{\frac{dt}{\alpha {{t}^{2}}-\beta }}=0 \qquad\text{ and }\qquad 
\frac 2\pi   \int\limits_{0}^{\infty }{\frac{dt}{\alpha {{t}^{2}}+\beta }}=\frac 1{  \sqrt{\alpha \beta }} .\\ 
\end{aligned}\]
Similarly, for function \[g\left( x \right)=\left\{ \begin{matrix}
   \sqrt{1-{{x}^{2}}} & \text{for }\left| x \right|<1  \\
   0 & \text{otherwise}  \\
\end{matrix} \right.\]
we have for $\left| x \right|<1$
\[   \tilde{g}\left( x \right)=\frac{1}{\pi }\int\limits_{-1}^{1}{\frac{\sqrt{1-{{y}^{2}}}}{x-y}}dy=\frac{1}{\pi }\int\limits_{-1}^{1}{\frac{1-{{y}^{2}}}{x-y}}\frac{dy}{\sqrt{1-{{y}^{2}}}} =\frac{1}{\pi }\int\limits_{-1}^{1}{\frac{{{x}^{2}}-{{y}^{2}}}{x-y}}\frac{dy}{\sqrt{1-{{y}^{2}}}} =x.
\]
 The real Hardy space ${{H}^{1}}\left( \mathbb{R} \right)$ is of all $f\in {{L}^{1}}\left( \mathbb{R} \right)$ such that $Hf\in {{L}^{1}}\left( \mathbb{R} \right).$ The duality of ${{H}^{1}}\left( \mathbb{R} \right)$ is $BMO\left( \mathbb{R} \right)$ the space of real functions of bounded mean oscillations \cite{Fefferman}. Clearly, ${{L}^{\infty }}\left( \mathbb{R} \right)\subseteq BMO\left( \mathbb{R} \right)$ but there are unbounded functions in  $BMO\left( \mathbb{R} \right)$ for example, the logarithmic function  $\ln x.$ We can define the logarithmic integral  
\[F\left( b \right)=\frac{1}{\pi }\int\limits_{-\infty }^{\infty }{f\left( x \right)\ln \frac{1}{\left| x-b \right|}}\cdot dx\]
 for a function $f\in {{H}^{1}}\left( \mathbb{R} \right)$ via ${{H}^{1}}-$ $BMO\left( \mathbb{R} \right)$  duality \cite{Fefferman}. Moreover, the Hilbert transform is a unitary operator acting on ${{L}^{2}}\left( \mathbb{R} \right).$ To compute the Hilbert transform of several functions  we define the complex Hardy sapces ${{\mathfrak{H}}^{p}}\left( {{\mathbb{C}}_{+}} \right)$ where ${{\mathbb{C}}_{+}}=\left\{ z\in \mathbb{C}:\text{ Im}\left( z \right)>0 \right\}$ and $1\le p\le \infty .$ More exactly, $\varphi \in {{\mathfrak{H}}^{p}}\left( {{\mathbb{C}}_{+}} \right)$ \cite{Duren}  if $\varphi $ is analytic in ${{\mathbb{C}}_{+}}$ and 
\[\left\| \varphi  \right\|_{p}^{p}:=\underset{y>0}{\mathop{\sup }}\,\int\limits_{-\infty }^{\infty }{{{\left| \varphi \left( x+iy \right) \right|}^{p}}}dx<\infty .\]
 If $p=\infty $ then ${{\mathfrak{H}}^{\infty }}\left( {{\mathbb{C}}_{+}} \right)$ is defined to be the space of bounded analytic function in ${{\mathbb{C}}_{+}}.$ It is well known that if $\varphi \in {{\mathfrak{H}}^{p}}\left( {{\mathbb{C}}_{+}} \right)$ then for almost every $x\in \mathbb{R}$ there is $\underset{y\to 0}{\mathop{\lim }}\,\varphi \left( x+iy \right)=:f\left( x \right)+i\tilde{f}\left( x \right),$ where $f,\tilde{f}\in {{L}^{p}}\left( \mathbb{R} \right)$ if $1<p<\infty .$ (Note that $f\left( x \right)=\operatorname{Re}\varphi \left( x+i0 \right)$ and $\tilde{f}\left( x \right)=\operatorname{Im}\varphi \left( x+i0 \right).$) Therefore, the Hilbert transform is bounded on ${{L}^{p}}\left( \mathbb{R} \right)$ for $1<p<\infty $ and $H\left( Hf \right)=-f$ for every $f\in {{L}^{p}}\left( \mathbb{R} \right)$ with $1<p<\infty .$  We have the formula
\[\varphi \left( z \right)=\frac{i}{\pi }\int\limits_{-\infty }^{\infty }{\frac{\operatorname{Re}\varphi \left( t \right)}{z-t}}\cdot dt=\frac{1}{\pi }\int\limits_{-\infty }^{\infty }{\frac{\operatorname{Im}\varphi \left( t \right)}{t-z}}\cdot dt\]
for any $\varphi \in {{\mathfrak{H}}^{p}}\left( {{\mathbb{C}}_{+}} \right)$ with  $1\le p<\infty .$ It is also known that $\tilde{f}$ is locally integrable if $f\in {{L}^{1}}\left( \mathbb{R} \right).$ On the other hand, we can define Hilbert transform of $f\in {{L}^{\infty }}\left( \mathbb{R} \right)$ up to a constant.  For example, $H\left( \cos x \right)=\sin x$ and $H\left( {{e}^{f}}\cos \tilde{f} \right)={{e}^{f}}\sin \tilde{f}$ for any $f\in {{L}^{\infty }}\left( \mathbb{R} \right).$ It is well known that if $f\in {{L}^{\infty }}\left( \mathbb{R} \right)$ then $\tilde{f}\in BMO\left( \mathbb{R} \right).$  Moreover, 
\[\int\limits_{-\infty }^{\infty }{f\left( x \right)\tilde{g}\left( x \right)dx}=-\int\limits_{-\infty }^{\infty }{\tilde{f}\left( x \right)g\left( x \right)dx}\]
  for $f\in {{L}^{p}}\left( \mathbb{R} \right)$ and $g\in {{L}^{q}}\left( \mathbb{R} \right)$ with $1<p<\infty $
and $\frac{1}{p}+\frac{1}{q}=1.$ Replace $g$ by ${{\chi }_{\left( a,b \right)}}$we have 
\[\frac{1}{\pi }\int\limits_{-\infty }^{\infty }{f\left( x \right)\ln \left| \frac{x-a}{x-b} \right|dx}=-\int\limits_{a}^{b}{\tilde{f}\left( x \right)dx}\]
 for every $f\in {{L}^{p}}\left( \mathbb{R} \right).$ For a rapidly decay function $f$ we can define  the logarithmic integral 
\[F\left( b \right)=\frac{1}{\pi }\int\limits_{-\infty }^{\infty }{f\left( x \right)\ln \frac{1}{\left| x-b \right|}}\cdot dx.\]
 Then 
\[F\left( b \right)-F\left( a \right)=-\int\limits_{a}^{b}{\tilde{f}\left( x \right)}dx.\]
 Hence, $F$ is locally absolutely continuous with weak derivative $-\tilde{f}.$ Moreover, if $\varphi \in {{\mathfrak{H}}^{2}}\left( {{\mathbb{C}}_{+}} \right)$ then ${{\varphi }^{2}}\in {{\mathfrak{H}}^{1}}\left( {{\mathbb{C}}_{+}} \right)$ and consequently, for almost every $x\in \mathbb{R},$ 
\[\underset{y\to 0}{\mathop{\lim }}\,\varphi {{\left( x+iy \right)}^{2}}={{\left[ f\left( x \right)+i\tilde{f}\left( x \right) \right]}^{2}}=f{{\left( x \right)}^{2}}-\tilde{f}{{\left( x \right)}^{2}}+2if\left( x \right)\tilde{f}\left( x \right).\]
 Thus, $H\left( {{f}^{2}}-{{{\tilde{f}}}^{2}} \right)=$ $2f\tilde{f}$ for every $f\in {{L}^{2}}\left( \mathbb{R} \right).$ Hence, ${{f}^{2}}-{{\tilde{f}}^{2}}$ and  $f\tilde{f}\in {{H}^{1}}\left( \mathbb{R} \right)$ for every $f\in {{L}^{2}}\left( \mathbb{R} \right).$ 
This is a tipical example for functions in ${H}^{1}\left( \mathbb{R} \right)$.
More generally, let $\varphi \in {{\mathfrak{H}}^{p}}\left( {{\mathbb{C}}_{+}} \right)$ and $\phi \in {{\mathfrak{H}}^{q}}\left( {{\mathbb{C}}_{+}} \right)$ with $\frac{1}{p}+\frac{1}{q}\le 1.$  Then $\varphi \phi \in {{\mathfrak{H}}^{r}}\left( {{\mathbb{C}}_{+}} \right)$ with $\frac{1}{r}=\frac{1}{p}+\frac{1}{q}$ so we have 
\[H\left( f\tilde{g}+\tilde{f}g \right)=\tilde{f}\tilde{g}-fg\quad\text{ with } f\in {{L}^{p}}\left( \mathbb{R} \right)\text{ and } g\in {{L}^{q}}\left( \mathbb{R} \right).\]
We define the Fourier transform $\widehat{f}$ of a function $f\in {{L}^{1}}\left( \mathbb{R} \right)$ by 
\[{\cal F}\left( f,x \right)=\widehat{f}\left( x \right)=\frac{1}{\sqrt{2\pi }}\int\limits_{-\infty }^{\infty }{f\left( t \right){{e}^{-itx}}dt,}\] for  $x\in \mathbb{R}.$
 Then  $\widehat{f}$ is uniformly continuous on  $\mathbb{R}$ and $\underset{\left| x \right|\to \infty }{\mathop{\lim }}\,\widehat{f}\left( x \right)=0.$ For $f\in {{H}^{1}}\left( \mathbb{R} \right),$ we have $\widehat{Hf}\left( x \right)=-i\widehat{f}\left( x \right)\cdot \text{sign }x$ for every $x\in \mathbb{R}$ and $\int\limits_{-\infty }^{\infty }{\left| \frac{\widehat{f}\left( x \right)}{x} \right|dx<\infty }$ (Hardy inequality).  Now consider the logarithmic integral 
\[F\left( b \right)=\frac{1}{\pi }\int\limits_{-\infty }^{\infty }{f\left( x \right)\ln \frac{1}{\left| x-b \right|}}\cdot dx\]
 of a function $f\in {{H}^{1}}\left( \mathbb{R} \right),$ which is defined via the duality ${{H}^{1}}-BMO\left( \mathbb{R} \right)$ \cite{Fefferman}.
We will prove that 
\[F\left( b \right)=-\int\limits_{-\infty }^{b}{\tilde{f}\left( x \right)dx}\]
 for every $b\in \mathbb{R}.$ It is enough to prove this equality for a rapidly decay function $f\in {{H}^{1}}\left( \mathbb{R} \right).$ As we have seen before, $F\left( b \right)-F\left( a \right)=-\int\limits_{a}^{b}{\tilde{f}\left( x \right)dx}.$ This means that $F\left( b \right)$ is locally absolutely continuous and ${F}'\left( b \right)=-\tilde{f}\left( b \right)$ so $\widehat{{{F}'}}=-\widehat{{\tilde{f}}}$ where the Fourier transform is taken in distributional sense. Consequently, $it\widehat{F}\left( t \right)=i\widehat{f}\left( t \right)\text{sign }t$ or equivalently, \[\widehat{F}\left( t \right)=\frac{\widehat{f}\left( t \right)}{\left| t \right|} \in {{L}^{1}}\left( \mathbb{R} \right)
\quad\hbox{(by Hardy inequality )}
\] and by inversion formula 
\[F\left( b \right)=\frac{1}{\pi }\int\limits_{-\infty }^{\infty }{f\left( x \right)\ln \frac{1}{\left| x-b \right|}}\cdot dx=\frac{1}{\sqrt{2\pi }}\int\limits_{-\infty }^{\infty }{\frac{\widehat{f}\left( t \right)}{\left| t \right|}\cdot {{e}^{ibt}}dt}\]
which is uniformly continuous on $\mathbb{R}$ and $\underset{\left| a \right|\to \infty }{\mathop{\lim }}\,F\left( a \right)=0$. Now from the formula 
\[F\left( b \right)-F\left( a \right)=-\int\limits_{a}^{b}{\tilde{f}\left( x \right)dx}\]
 we have 

\bigskip
\par\noindent {\bf Theorem 1.}  {\it For any function $f\in {{H}^{1}}\left( \mathbb{R} \right)$ and $b\in \mathbb{R},$
\[\frac{1}{\pi }\int\limits_{-\infty }^{\infty }{f\left( x \right)\ln \frac{1}{\left| x-b \right|}}\cdot dx=-\int\limits_{-\infty }^{b}{\tilde{f}\left( x \right)dx}.\]   }

\bigskip
\par\noindent{\bf Remark. }
It is proved in \cite{Stefanov}  that if $f\in {{H}^{1}}\left( \mathbb{R} \right)$ then the logarithmic integral $F$  is of bounded variation. Our result is much stronger.  Now replace $f\in {{H}^{1}}\left( \mathbb{R} \right)$ by ${{f}^{2}}-{{\tilde{f}}^{2}}$ we have

\bigskip
\par\noindent {\bf Theorem 2.}  {\it For any function $f\in {{L}^{2}}\left( \mathbb{R} \right)$ and $b\in \mathbb{R},$
$$\frac{1}{\pi }\int\limits_{-\infty }^{\infty }{\left[ {{f}^{2}}\left( x \right)-{{{\tilde{f}}}^{2}}\left( x \right) \right]\ln \frac{1}{\left| x-b \right|}}\cdot dx=-2\int\limits_{-\infty }^{b}{f\left( x \right)\tilde{f}\left( x \right)dx}.$$  }

\bigskip
\par\noindent 
For example, take 
\[\varphi \left( z \right)=\frac{i}{z+i}\in {{\mathfrak{H}}^{2}}\left( {{\mathbb{C}}_{+}} \right)
\]
 then 
\[\varphi \left( x \right)=\frac{i}{x+i}=\frac{i\left( x-i \right)}{{{x}^{2}}+1}=\frac{1}{{{x}^{2}}+1}+i\cdot \frac{x}{{{x}^{2}}+1}\]
   so
\[f\left( x \right)=\frac{1}{{{x}^{2}}+1}
\quad\hbox{ and }\quad \tilde{f}\left( x \right)=\frac{x}{{{x}^{2}}+1}\in {{L}^{2}}\left( \mathbb{R} \right).\]
 Therefore, 
\[\frac{1-{{x}^{2}}}{{{\left( {{x}^{2}}+1 \right)}^{2}}},
\quad
\frac{x}{{{\left( {{x}^{2}}+1 \right)}^{2}}}\in {{H}^{1}}\left( \mathbb{R} \right)\]
 and 
\[H\left( \frac{1-{{x}^{2}}}{{{\left( {{x}^{2}}+1 \right)}^{2}}} \right)=\frac{2x}{{{\left( {{x}^{2}}+1 \right)}^{2}}}=-{f}'\left( x \right).\]
 Hence,
\[\frac{1}{\pi }\int\limits_{-\infty }^{\infty }{\frac{1-{{x}^{2}}}{{{\left( {{x}^{2}}+1 \right)}^{2}}}\ln \frac{1}{\left| x-b \right|}}\cdot dx=\int\limits_{-\infty }^{b}{{f}'\left( x \right)dx}=f\left( b \right)=\frac{1}{{{b}^{2}}+1}.\]
Clearly, $\tilde{f}\left( x \right)=\frac{x}{{{x}^{2}}+1}\notin {{L}^{1}}\left( \mathbb{R} \right)$
 so $f\left( x \right)=\frac{1}{{{x}^{2}}+1}\notin {{H}^{1}}\left( \mathbb{R} \right).$  Now note 
that the function $\varphi \left( z \right)=\frac1{\sqrt{1-z^2}}$ is in 
${{\mathfrak{H}}^{p}}\left( {{\mathbb{C}}_{+}} \right)$ 
for any $p\in \left( 1,2 \right)$ but it does not belong to 
${{\mathfrak{H}}^{1}}\left( {{\mathbb{C}}_{+}} \right)\cup {{\mathfrak{H}}^{2}}\left( {{\mathbb{C}}_{+}} \right).$ 
(The square root is taken in the sense that the real part of $\varphi $ is positive.) Indeed, we have
 \[
\varphi \left( x+i0 \right)=\frac{1}{\sqrt{1-{{x}^{2}}}}=f\left( x \right)+i\tilde{f}\left( x \right)
\]
 with $f\tilde{f}=0$ because $f\left( x \right)=\frac{1}{\sqrt{1-{{x}^{2}}}}$ for 
$\left| x \right|<1$ and $f\left( x \right)=0$  for $\left| x \right|>1.$  Similarly, 
$\tilde{f}\left( x \right)=0$ for $\left| x \right|<1$ and $\tilde{f}\left( x \right)=\frac{1}{\sqrt{{{x}^{2}}-1}}$ for $x>1$ and $\tilde{f}\left( x \right)=-\frac{1}{\sqrt{{{x}^{2}}-1}}$  for $x<-1.$  Thus, $f\tilde{f}=0$ and ${{f}}\notin {{L}^{2}}\left( \mathbb{R} \right).$ 
Now we note that $\varphi \left( z \right)={{e}^{-{{z}^{2}}}}$  is analytic on the complex plane $\mathbb{C}$ but it does not belong to any ${{\mathfrak{H}}^{p}}\left( {{\mathbb{C}}_{+}} \right).$ Indeed, if otherwise the boundary function 
\[\varphi \left( x+i0 \right)={{e}^{-{{x}^{2}}}}=f\left( x \right)+i\tilde{f}\left( x \right)\] 
with $f\left( x \right)={{e}^{-{{x}^{2}}}}$ and $\tilde{f}\left( x \right)=0$ which is absurd. On the other hand, \[\int\limits_{-\infty }^{\infty }{{{\left| \varphi \left( x+iy \right) \right|}^{p}}}dx={{e}^{p{{y}^{2}}}}\int\limits_{-\infty }^{\infty }{{{e}^{-p{{x}^{2}}}}}dx\to \infty \]
 as $y\to \infty ,$  which means that $\varphi \notin {{\mathfrak{H}}^{p}}\left( {{\mathbb{C}}_{+}} \right).$  
 Finally, let \[{{a}_{1}}<{{a}_{2}}<\cdots <{{a}_{2\ell }},
\quad
E=\bigcup\limits_{k=1}^{\ell }{\left[ {{a}_{2k-1}},{{a}_{2k}} \right]}\text{  and }
K\left( x \right)=\prod\limits_{j=1}^{2\ell }{\left( x-{{a}_{j}} \right)}.\] 
Then $K(x)\le 0$ if and only if $x\in E.$ Let
\[g\left( x \right)=g_E(x)=\left\{ \begin{matrix}   (-1)^{\ell-k}\sqrt{\left| K\left( x \right) \right|} & \quad\text{ if }x\in  {\left[ {{a}_{2k-1}},{{a}_{2k}} \right]} \\   0 & \quad\text{otherwise}\text{.}  \\\end{matrix} \right.\]
Note that 
\[\varphi \left( z \right)=\frac{{{z}^{k-1}}}{\sqrt{K\left( z \right)}}\in {{\mathfrak{H}}^{p}}\left( {{\mathbb{C}}_{+}} \right)\text{  for  }k=1,2,\cdots ,\ell \text{ and }p\in \left( 1,2 \right).\]
Here, $\sqrt{K\left( z \right)}\sim {{z}^{\ell }}$ as $z\to \infty .$ 
Moreover,
\[\varphi \left( x \right)=\left\{ \begin{matrix}   (\frac{{(-1)^{\ell-m}{x}^{k-1}}}{\sqrt{K(x)}} & \text{ if }& x\in  (a_{2m},a_{2m+1}) \\   -\frac{i{{x}^{k-1}}}{g\left( x \right)}  &\text{if} & x\in E.\\
\end{matrix} \right.\]
($m=0,1,\cdots,\ell$, $a_0=-\infty$ and $a_{2\ell+1}=\infty$). In fact, it is easily follows from the computation of the positive harmonic argument of $K(z)$. More exactly, Arg$K(x)=(2\ell-j)\pi$ for $x\in(a_j,a_{j+1})$ and $j=0,1,\cdots,2\ell$.
Thus,
 \[\frac1\pi \int_E \frac {y^{k-1}}{g_E(y)}\frac{dy}{x-y}=\left\{ \begin{matrix}   (\frac{{(-1)^{\ell-m}{x}^{k-1}}}{\sqrt{K(x)}} & \text{ if }& x\in  (a_{2m},a_{2m+1}) \\   0 & \text{ if } & x\in E.  \\\end{matrix} \right.\]
For example, 
\[\underset{\operatorname{Im}z\to 0}{\mathop{\lim }}\,\frac{1}{\sqrt{\left( z-a \right)\left( z-b \right)}}=\left\{ -\begin{matrix}
   \frac{1}{\sqrt{\left( x-a \right)\left( x-b \right)}} & \text{if} & x>b  \\
   \frac{i}{\sqrt{\left( x-a \right)\left( b-x \right)}} & \text{if} & a<x<b  \\
   -\frac{1}{\sqrt{\left( x-a \right)\left( x-b \right)}} & \text{if} & x<a  \\
\end{matrix} \right.
\quad \left( x=\operatorname{Re}z \right)
\]
so we have
\[\frac{1}{\pi }\int\limits_{a}^{b}{\frac{1}{\sqrt{\left( y-a \right)\left( b-y \right)}}\frac{dy}{x-y}}=\left\{ \begin{matrix}
   \frac{1}{\sqrt{\left( x-a \right)\left( x-b \right)}} & \text{if} & x>b  \\
   0 & \text{if} & a<x<b  \\
   -\frac{1}{\sqrt{\left( x-a \right)\left( x-b \right)}} & \text{if} & x<a.  \\
\end{matrix} \right.\]
Similarly, 
\[\underset{\operatorname{Im}z\to 0}{\mathop{\lim }}\,z-\sqrt{\left( z-a \right)\left( z-b \right)}=\left\{ \begin{matrix}
   x-\sqrt{\left( x-a \right)\left( x-b \right)} & \text{if} & x>b  \\
   x+\sqrt{\left( x-a \right)\left( x-b \right)} & \text{if} & x<a  \\
   x-i\sqrt{\left( x-a \right)\left( b-x \right)} & \text{if} & a<x<b  \\
\end{matrix} \right.\text{    }\left( x=\operatorname{Re}z \right)\]
and we have
\[\frac{1}{\pi }\int\limits_{a}^{b}{\frac{\sqrt{\left( y-a \right)\left( b-y \right)}}{x-y}}dy=\left\{ \begin{matrix}
   x-\frac{a+b}2-\sqrt{\left( x-a \right)\left( x-b \right)} & \text{if} & x>b  \\
   x-\frac{a+b}2+\sqrt{\left( x-a \right)\left( x-b \right)} & \text{if} & x<a  \\
   x -\frac{a+b}2& \text{if} & a<x<b.  \\
\end{matrix} \right.\]
More generally, let $E=\left[ -b,-a \right]\cup \left[ a,b \right]$ and
\[g\left( x \right)=\left\{ \begin{matrix}
   -\sqrt{\left( {{b}^{2}}-{{x}^{2}} \right)\left( {{x}^{2}}-{{a}^{2}} \right)} & \text{if} & x\in \left[ -b,-a \right]  \\
   \sqrt{\left( {{b}^{2}}-{{x}^{2}} \right)\left( {{x}^{2}}-{{a}^{2}} \right)} & \text{if} & x\in \left[ a,b \right]  \\
   0 & \text{otherwise}\text{.} & {}  \\
\end{matrix} \right.\]
Then  
\[\begin{aligned}
 \tilde{g}\left( x \right)&=\frac{1}{\pi }\int\limits_{a}^{b}{\sqrt{\left( {{b}^{2}}-{{y}^{2}} \right)\left( {{y}^{2}}-{{a}^{2}} \right)}\left[ \frac{1}{x-y}-\frac{1}{x+y} \right]dy} \\ 
 & =\frac{1}{\pi }\int\limits_{a}^{b}{\sqrt{\left( {{b}^{2}}-{{y}^{2}} \right)\left( {{y}^{2}}-{{a}^{2}} \right)}\frac{2ydy}{\left( {{x}^{2}}-{{y}^{2}} \right)}} \\ 
 & =\frac{1}{\pi }\int\limits_{{{a}^{2}}}^{{{b}^{2}}}\frac{\sqrt{\left( {{b}^{2}}-y \right)\left( y-{{a}^{2}} \right)}}{\left( {{x}^{2}}-y \right)}dy\\
&=\left\{ \begin{matrix}
   x^2-\frac{a^2+b^2}2-\sqrt{\left( x^2-a^2 \right)\left( x^2-b^2 \right)} & \text{if} & |x|>b  \\
   x^2-\frac{a^2+b^2}2+\sqrt{\left( x^2-a^2 \right)\left( x^2-b^2 \right)} & \text{if} &|x| <a  \\
   x^2 -\frac{a^2+b^2}2& \text{if} & a<|x|<b.  \\
\end{matrix} \right.
\end{aligned}
\]

\bigskip
\bigskip

\sect{\bf Finite Hilbert transforms and Inversion}

\bigskip
\par\noindent Now we are interested in compactly supported positive functions and their Hilbert transforms. More exactly, if $f$ is supported in the interval $[a,b]$, the finite Hilbert transform of $f$ is given by the Cauchy principal value integral $Hf(s)=(1/\pi)\int_a^b(s-t)^{-1}f(t)\,dt$ for real $s$. By complex variable arguments we have the inversion formula  \cite{Tricomi} \cite{You}  
\[f(t)=\frac{1}{\pi \sqrt{(t-a)(b-t)}}\left( \int_{a}^{b}{\frac{Hf(s)}{s-t}}\sqrt{(s-a)(b-s)}ds+\int_{a}^{b}{f}(s)ds \right)\]
for $f\in L^p(\R)$  with $p>1$. 
Now let 
	\[E=\bigcup\limits_{k=1}^{\ell }{\left[ {{a}_{2k-1}},{{a}_{2k}} \right]}\]
 be the finite  union of intervals and assume that $f$ is supported in  $E$.  We are interested in the inversion formula of the Hilbert transform of $f$. To this end, let 
	\[K\left( x \right)=\prod\limits_{j=1}^{2\ell }{\left( x-{{a}_{j}} \right)}\quad\text{and }\]
\[g\left( x \right)=g_E(x)=\left\{ \begin{matrix}
   (-1)^{\ell-k}\sqrt{\left| K\left( x \right) \right|} & \quad\text{ if }x\in  {\left[ {{a}_{2k-1}},{{a}_{2k}} \right]} \\
   0 & \quad\text{otherwise}\text{.}  \\
\end{matrix} \right.
\]
 Then $K(x)\le 0$ for $x\in E$. 
Using the formula
\[H\left( f\tilde{g}+\tilde{f}g \right)=\tilde{f}\tilde{g}-fg\quad\text{ with } f\in {{L}^{p}}\left( \mathbb{R} \right)\text{ and } g\in {{L}^{q}}\left( \mathbb{R} \right)\]
  we have 
\[\begin{aligned}
 H\left( f\tilde{g}+g\tilde{f},x \right) 
 &=\frac{1}{\pi }\int\limits_{E}{\frac{\tilde{g}\left( y \right)f\left( y \right)}{x-y}}dy+\frac{1}{\pi }\int\limits_{E}{\frac{g\left( y \right)\tilde{f}\left( y \right)}{x-y}}dy \\ 
 & =-\frac{1}{\pi }\int\limits_{E}\frac{\tilde{g}\left( x \right)-\tilde{g}\left( y \right)}{x-y} f\left( y \right)dy+\tilde{g}\left( x \right)\tilde{f}\left( x \right)+\frac{1}{\pi }\int\limits_{E}\frac{g(y)\tilde{f}\left( y \right)}{x-y}dy \\ 
 & =\tilde{g}\left( x \right)\tilde{f}\left( x \right)-f\left( x \right){g\left( x \right) }\quad\text{  if }x\in E.
\end{aligned}\]
Therefore, 
the inversion formula
\[f\left( x \right)=\frac{1}{\pi g_E\left( x \right) } \left( \int\limits_{E}\frac{\tilde{g}_E\left( x \right)-\tilde{g}_E\left( y \right)}{x-y} f\left( y \right)dy+\int\limits_{E}\frac{g_E\left( y \right)\tilde{f}\left( y \right)}{y-x}dy \right)\]
holds for $f\in L^p$ with $p>1$ and $x\in E$. But this formula is not useful, because we should compute the Hilbert transform $\tilde g$.  In the next section we will prove that  that $\tilde g_E$ on $E$ is a polynomial of degree $\ell$. Hence, the first term
\[\int\limits_{E}\frac{\tilde{g}_E\left( x \right)-\tilde{g}_E\left( y \right)}{x-y} f\left( y \right)dy\]
is a polynomial of degree $\le\ell-1$ which is determined uniquely by the first $\ell$  moments of $f$. 
It follows at once from this formula that if $\tilde f=0$ on $E$ then $f$ has the form $f(x)=\rho(x)/g_E(x)$, where $\rho$ is a polynomial of degree less than $\ell$. On the other hand, at the end of section 1, we have seen that every function of this form (supported in the set $E$) has Hilbert transform vanishing in $E$.
We make another inversion formula which is more applicable. Recall that the equilibrium measure of  a compact set $E$ is the only solution of the energy optimization problem 
\[I\left( \mu  \right)=\iint{\ln \frac{1}{\left| x-t \right|}}d\mu \left( x \right)d\mu \left( t \right)\to \min \]
subject to every Borel probability measure $\mu$ supported in $E$. 
The density function $\omega_E $ of the  equilibrium measure of $E$ is 
	\[{{\omega }_{E}}(x)=\frac{1}{\pi }\cdot \frac{\left| {{\rho }_{\ell -1}}\left( x \right) \right|}{\sqrt{\left| K\left( x \right) \right|}}=\frac{1}{\pi }\cdot \frac{{{\rho }_{\ell -1}}\left( x \right)}{g\left( x \right)}\]
 where ${{\rho }_{\ell -1}}\left( x \right)={{x}^{\ell -1}}+\cdots=(t-\tau_1 )(t-\tau_2 )\cdots(t-\tau_{\ell -1} )$ is that unique polynomial satisfying 
	\[\int\limits_{{{a}_{2j}}}^{{{a}_{2j+1}}}{\frac{{{\rho }_{\ell -1}}\left( x \right)}{\sqrt{\left| K\left( x \right) \right|}}}\cdot dx=0\]
 for $j=1,2,\cdots ,\ell -1$ and 
\[g\left( x \right)=\left\{ \begin{matrix}
   (-1)^{\ell-k}\sqrt{\left| K\left( x \right) \right|} & \quad\text{ if }x\in  {\left[ {{a}_{2k-1}},{{a}_{2k}} \right]} \\
   0 & \quad\text{otherwise}\text{.}  \\
\end{matrix} \right.\] 
The roots of ${{\rho }_{\ell -1}}$ are in the gaps of $E$. 
 For example, 
\[{{\omega }_{\left[ a,b \right]}}\left( x \right)=\frac{1}{\pi \sqrt{\left( x-a \right)\left( b-x \right)}}\]
and
\[{{\omega }_{\left[ -b,-a \right]\cup \left[ a,b \right]}}\left( x \right)=\frac{\left| x \right|}{\pi \sqrt{\left( {{x}^{2}}-{{a}^{2}} \right)\left( {{b}^{2}}-{{x}^{2}} \right)}}.\]
Moreover, it is well known that the Hilbert transform of the density function $\omega_E$
is zero in $E$. It follows from the fact that
\[\int_E \frac {y^{k-1}}{g(y)}\frac{dy}{x-y}=0\quad\text{ for } x\in E
\text{ and } k=1,2,\cdots,\ell. 
\]
The density function itself is in $L^q$ for any $q<2$.
Let $g_0=\omega_E$ and we try to use the formula $H\left( f\tilde{g_0}+\tilde{f}g_0 \right)=\tilde{f}\tilde{g_0}-fg_0$. Because $g_0\in {{L}^{q}}$ for any $q<2$ we should assume that $f\in {{L}^{p}}$ with $p>2.$ On the other hand, $f\tilde{g_0}$ is identically 0, because $f$ is supported on $E$ and $\tilde{g_0}=0$ on $E.$ Hence, $H\left( \tilde{f}g_0,x \right)=-f\left( x \right)g_0\left( x \right)$ for $x\in E.$ Therefore, we have

\bigskip
\par\noindent {\bf Theorem 3.}  {\it
Let $f\in {{L}^{p}}\left( \mathbb{R} \right)$ for some $p> 2.$ If $f$ is supported in a compact set $E$ then
\[f\left( x \right)=\frac{1}{\pi {{\omega }_{E}}\left( x \right)}\int\limits_{E}{\frac{\tilde{f}\left( y \right){{\omega }_{E}}\left( y \right)}{y-x}}dy\quad\text{     for a.e. } x\in E,\]
where $\omega_E$ denotes the density function of the equilibrium measure of $E.$ The compact set $E$ is supposed to have absolutely continuous equilibrium measure and the density function $\omega_E$ is supposed to be in ${{L}^{q}}\left( \mathbb{R} \right)$ with  $\frac 1p+ 
\frac 1q=1$.
}

\bigskip
\par\noindent {\bf Remark. } The assumption $p>2$ is very essential. Otherwise, the density function $\omega_E$ itself does not satisfy this inversion formula. 
Moreover, if we take $g_0(x)=\rho(x)/g(x)$ for $x\in E$ and $g_0(x)=0$ for $x\notin E$ we also have $\tilde g_0(x)=0$ for $x\in E$. Here, $\rho$ denotes a polynomial of degree $<\ell$ (the number of holes of $E$). 
Therefore, if $f $ is supported in $E$ and $f\in L^p$ for some $p>2$ then
 \[f\left( x \right)=\frac{g(x)}{\pi \rho(x)}\int\limits_{E}\frac{\rho(y)\tilde{f}\left( y \right)}{g(y)(y-x)}dy\quad\text{     for a.e. } x\in E.\]
Recall that
	\[K\left( x \right)=\prod\limits_{j=1}^{2\ell }{\left( x-{{a}_{j}} \right)}\quad\text{and }
g\left( x \right)=\left\{ \begin{matrix}
   (-1)^{\ell-k}\sqrt{\left| K\left( x \right) \right|} & \quad\text{ if }x\in  {\left[ {{a}_{2k-1}},{{a}_{2k}} \right]} \\
   0 & \quad\text{otherwise}\text{.}  \\
\end{matrix} \right.
\] 

\bigskip
\noindent
For example, let $E=[-1,1]$ and
\[g\left( x \right)=\left\{ \begin{matrix}
   \sqrt{1-{{x}^{2}}} & \quad\text{if }\left| x \right|<1  \\
   0 & \quad\text{otherwise}  \\
\end{matrix} \right.\]
Then $\tilde{g}\left( x \right)=x$  if  $\left| x \right|<1.$   On the other hand,
\[H\left[ {{\left( \tilde{g}+ig \right)}^{n}} \right]=i{{\left( \tilde{g}+ig \right)}^{n}}  \quad \text{ and }\quad H\left[ {{\left( \tilde{g}-ig \right)}^{n}} \right]=-i{{\left( \tilde{g}-ig \right)}^{n}}\]
so 
\[\begin{aligned}
H\left[ {{\left( \tilde{g}+ig \right)}^{n}}-{{\left( \tilde{g}-ig \right)}^{n}} \right]  & =i\left[ {{\left( \tilde{g}+ig \right)}^{n}}+{{\left( \tilde{g}-ig \right)}^{n}} \right] \\ 
 & =2i{{T}_{n}}\left( x \right)\quad\text{  for }\left| x \right|<1  
\end{aligned}\]
where ${{T}_{n}}\left( x \right)=\cos n\theta $ is Chebisev polynomial of first kind. Moreover,
\[\begin{aligned}
 {{\left( \tilde{g}+ig \right)}^{n}}-{{\left( \tilde{g}-ig \right)}^{n}} 
& =2ig\left[ {{\left( \tilde{g}+ig \right)}^{n-1}}+\cdots +{{\left( \tilde{g}-ig \right)}^{n}} \right] \\ 
 & =0\quad\text{  if  }\left| x \right|>1 \\ 
 & =2i{{U}_{n-1}}\left( x \right)\sqrt{1-{{x}^{2}}}\quad\text{  for }\left| x \right|<1,  
\end{aligned}\]
where \[{{U}_{n-1}}\left( x \right)=\frac{\sin n\theta }{\sin \theta }\] is Chebisev polynomial of second kind. Therefore, 
\[\frac{1}{\pi }\int\limits_{-1}^{1}{\frac{{{U}_{n-1}}\left( x \right)\sqrt{1-{{x}^{2}}}}{x-y}}dx=-{{T}_{n}}\left( y \right)\quad\text{  for }\left| y \right|<1,\]
and in virtue of inversion of finite Hilbert transform
\[\frac{1}{\pi }\int\limits_{-1}^{1}{\frac{{{T}_{n}}\left( x \right)}{x-y}}\frac{dx}{\sqrt{1-{{x}^{2}}}}={{U}_{n-1}}\left( y \right)\quad\text{  for }\left| y \right|<1.\]
We get at least
\[\frac{1}{\pi }\int\limits_{-1}^{1}{\frac{{{T}_{n}}\left( x \right)-{{T}_{n}}\left( y \right)}{x-y}}\frac{dx}{\sqrt{1-{{x}^{2}}}}={{U}_{n-1}}\left( y \right)\quad\text{  for }\left| y \right|<1\]
and
\[\frac{1}{\pi }\int\limits_{-1}^{1}{\frac{{{U}_{n}}\left( x \right)-{{U}_{n}}\left( y \right)}{x-y}}\sqrt{1-{{x}^{2}}}dx={{U}_{n-1}}\left( y \right)\quad\text{  for }\left| y \right|<1.\]
Now let $E=\left[ -1,-\alpha  \right]\cup \left[ \alpha ,1 \right]$  where
$\alpha $ maximizes the function
\[\Phi \left( \tau  \right)=\frac12\ln \frac{1-\tau^2 }{4}+\frac{a}{\pi }\int\limits_{\tau }^{1}{\frac{{{x}^{3}}dx}{\sqrt{\left( 1-{{x}^{2}} \right)\left( {{x}^{2}}-{{\tau }^{2}} \right)}}}
=\frac12\ln \frac{1-\tau^2 }{4}+\frac{a(\tau^2+1)}4
\]
and $a>2$ is fixed. Then $\alpha=\sqrt{\frac{a-2}a}$. 
 Suppose that $f$ is supported in $E$    and $\tilde{f}\left( x \right)=-ax$ for $x\in E$. 
Let
\[g\left( x \right)=\left\{ \begin{matrix}
   -\sqrt{\left( 1-{{x}^{2}} \right)\left( {{x}^{2}}-{{\alpha }^{2}} \right)} & \text{if} & x\in \left[ -1,-\alpha  \right]  \\
   \sqrt{\left( 1-{{x}^{2}} \right)\left( {{x}^{2}}-{{\alpha }^{2}} \right)} & \text{if} & x\in \left[ \alpha ,1 \right]  \\
   0 & \text{otherwise}\text{.} & {}  \\
\end{matrix} \right.\]
Then $\tilde{g}\left( x \right)={{x}^{2}}-\frac{\alpha^2+1}2$  for $x\in E$. Moreover, if $f\in L^p$ with $p>1$ and
 \[\frac{1}{\pi }\int\limits_{E}{f\left( x \right)dx=1}
\]
  then
\[f\left( x \right)=\frac{1}{\pi g\left( x \right)}\left( \pi x+\int\limits_{E}{\frac{ayg\left( y \right)}{x-y}}dy \right)=\sqrt{a}|x|\sqrt{\frac{ax^2+2-a}{1-x^2}}
\]
is positive in $E$.
\bigskip
\par\noindent {\bf Theorem 4.}  {\it
Let $f\in {{L}^{p}}\left( \mathbb{R} \right)$ for some $p> 4.$ If $f$ is supported in a compact set $E$ then
\[\int\limits_{E}{{{\left| f\left( y \right) \right|}^{2}}{{\omega }_{E}}\left( y \right)}dy=\int\limits_{E}{{{\left| \tilde{f}\left( y \right) \right|}^{2}}{{\omega }_{E}}\left( y \right)}dy\]
\[\int\limits_{E}{{{\left| f\left( y \right) \right|}^{2}}{{\omega }_{E}}\left( y \right)}ydy=\int\limits_{E}{{{\left| \tilde{f}\left( y \right) \right|}^{2}}{{\omega }_{E}}\left( y \right)y}dy\]
where $\omega_E$ denotes the density function of the equilibrium measure of $E.$ The compact set $E$ is supposed to have absolutely continuous equilibrium measure and the density function $\omega_E$ is supposed to be in ${{L}^{q}}\left( \mathbb{R} \right)$ with  $\frac 2p+\frac 1q=1$.
}

\bigskip
\par\noindent {\sl Proof: } Without loss of generality we assume that $f$  is real valued. Then
\[\begin{aligned}
 \int\limits_{E}{\left[ {{\left| \tilde{f}\left( y \right) \right|}^{2}}-{{\left| f\left( y \right) \right|}^{2}} \right]{{\omega }_{E}}\left( y \right)}dy & =\int\limits_{E}{\left[ \tilde{f}{{\left( y \right)}^{2}}-f{{\left( y \right)}^{2}} \right]{{\omega }_{E}}\left( y \right)}dy \\ 
 & =2\int\limits_{E}{H\left( \tilde{f}f,y \right){{\omega }_{E}}\left( y \right)}dy \\ 
 & =-2\int\limits_{\mathbb{R}}{\tilde{f}\left( y \right)f\left( y \right)H{{\omega }_{E}}\left( y \right)}dy=0  
\end{aligned}\]
(the Hilbert transform of ${{\omega }_{E}}$ is identically 0 on $E$) and the first identity is proved. For the second one, note that the Hilbert transform of $x{{\omega }_{E}}\left( x \right)$ is identically $-\frac1\pi$ on $E$ so
\[\begin{aligned}
  \int\limits_{E}{\left[ {{\left| \tilde{f}\left( y \right) \right|}^{2}}-{{\left| f\left( y \right) \right|}^{2}} \right]{{\omega }_{E}}\left( y \right)}ydy &=\int\limits_{E}{\left[ \tilde{f}{{\left( y \right)}^{2}}-f{{\left( y \right)}^{2}} \right]{{\omega }_{E}}\left( y \right)}ydy \\ 
 & =2\int\limits_{E}{H\left( \tilde{f}f,y \right)y{{\omega }_{E}}\left( y \right)}dy \\ 
 & =-2\int\limits_{\mathbb{R}}{\tilde{f}\left( y \right)f\left( y \right)H\left( y{{\omega }_{E}}\left( y \right) \right)}dy \\ 
 & =\frac2\pi\int\limits_{\mathbb{R}}{\tilde{f}\left( y \right)f\left( y \right)}dy=0.  
\end{aligned}\]
The proof is now complete.

\bigskip
\par\noindent {\bf Remark. } This theorem is proved in \cite{Rovnyak} in very special case where $f$ is continuous and $E=[-b,b]$. 

\bigskip
\bigskip

\sect{ \bf Orthonormal Polynomials and 
Finite Hilbert Transforms}

\bigskip
\noindent
Let $\{p_0,p_1,p_2,\cdots\}$ be the system of orthonormal polynomials with respect to the equilibrium measure of 
\[E=\bigcup\limits_{k=1}^{\ell }{\left[ {{a}_{2k-1}},{{a}_{2k}} \right]}.\]
Recall that 
\[K\left( x \right)=\prod\limits_{j=1}^{2\ell }{\left( x-{{a}_{j}} \right)},\]
\[g\left( x \right)=g_E(x)=\left\{ \begin{matrix}
   (-1)^{\ell-k}\sqrt{\left| K\left( x \right) \right|} & \quad\text{ if }x\in  {\left[ {{a}_{2k-1}},{{a}_{2k}} \right]} \\
   0 & \quad\text{otherwise,}  \\
\end{matrix} \right.\] 
and \[{{\omega }_{E}}(x)=\frac{1}{\pi }\cdot \frac{\left| {{\rho }_{\ell -1}}\left( x \right) \right|}{\sqrt{\left| K\left( x \right) \right|}}=\frac{1}{\pi }\cdot \frac{{{\rho }_{\ell -1}}\left( x \right)}{g\left( x \right)}\]
 where ${{\rho }_{\ell -1}}\left( x \right)={{x}^{\ell -1}}+\cdots=(t-\tau_1 )(t-\tau_2 )\cdots(t-\tau_{\ell -1} )$ is that unique polynomial satisfying 
	\[\int\limits_{{{a}_{2j}}}^{{{a}_{2j+1}}}{\frac{{{\rho }_{\ell -1}}\left( x \right)}{\sqrt{\left| K\left( x \right) \right|}}}\cdot dx=0\]
 for $j=1,2,\cdots ,\ell -1$.
Then we have a linear recurrence for $p_n$'s
\[x{{p}_{n}}(x)={{\alpha }_{n-1}}{{p}_{n-1}}(x)+{{\beta }_{n}}{{p}_{n}}(x)+{{\alpha }_{n}}{{p}_{n+1}}(x)\qquad \text{ for    }n=1,2,\cdots .\]
Here, $\left\{ {{\alpha }_{0}},{{\alpha }_{1}},\cdots  \right\}$ is a  bounded positive sequence and $\left\{ {{\beta }_{0}},{{\beta }_{1}},\cdots  \right\}$ is a  bounded real sequence. 
Moreover, ${{p}_{0}}=1$  and  
\[{{p}_{1}}\left( x \right)=\frac{x-{{\beta }_{0}}}{{{\alpha }_{0}}}.\]
Let
\[\begin{aligned}
 {{q}_{n-1}}\left( y \right) & =\int\limits_{E}{\frac{{{p}_{n}}\left( x \right)-{{p}_{n}}\left( y \right)}{x-y}{{\omega }_{E}}\left( x \right)dx}\quad\text{   for     }n=1,2,\cdots . \\ 
 \end{aligned}\]
Then $q_0=1/\alpha_0$, 
\[ {{q}_{n-1}}\left( y \right) 
  =-H\left( {{p}_{n}}{{\pi\omega }_{E}},y \right)\quad\text{    for  }  y\in E,   \quad n=1,2,\cdots\]
 and 
 \[\frac{{{q}_{n-1}}\left( x \right)}{{{p}_{n}}\left( x \right)}\to  \pi{{\tilde{\omega }}_{E}}\left( x \right)\quad\text{   as } n\to \infty \quad\text {  for  } x\notin E.\] 
Moreover, we have the linear recurrence for $q_n$'s (shifted one step)
\[x{{q}_{n}}(x)={{\alpha }_{n}}{{q}_{n-1}}(x)+{{\beta }_{n+1}}{{q}_{n}}(x)+{{\alpha }_{n+1}}{{q}_{n+1}}(x)\qquad \text{ for     } n=1,2,\cdots .\]
From the Christoffel-Darboux formula
\[\sum\limits_{k=0}^{n}{{{p}_{k}}\left( x \right){{p}_{k}}\left( y \right)}={{\alpha }_{n}}\frac{{{p}_{n+1}}\left( x \right){{p}_{n}}\left( y \right)-{{p}_{n}}\left( x \right){{p}_{n+1}}\left( y \right)}{x-y}\]
we have
\[\int\limits_{E}{\frac{{{p}_{n+1}}\left( x \right){{p}_{n}}\left( y \right)-{{p}_{n}}\left( x \right){{p}_{n+1}}\left( y \right)}{x-y}{{\omega }_{E}}\left( x \right)dx}=\frac{1}{{{\alpha }_{n}}}\]
and consequently,
\[\begin{aligned}
 p_{n+1} H\left( p_n \pi\omega_E \right) & =p_n H\left( p_{n+1}\pi \omega_E \right)+\frac 1{\alpha_n},  \\ 
 p_{n+1} q_{n-1}& =p_n q_n-\frac 1{\alpha_n}.\\ 
\end{aligned}\]
To compute the Hilbert transform of $g$ note that
\[\int_E \frac {y^{k-1}}{g(y)}\frac{dy}{x-y}=0\quad\text{ for } x\in E
\text{ and } k=1,2,\cdots,\ell,\]
and
\[g(y)^2=\left[ \sum\limits_{j=0}^{\ell +1}{{{\xi }_{j}}{{p}_{j}}\left( y \right)} \right]
{\rho }_{\ell -1}\left( y \right)+r(y),\quad\text{deg}(r)<\ell-1,
\]
so
\[\begin{aligned}
 \tilde{g}\left( x \right)  &=\frac{1}{\pi }\int\limits_{E}{\frac{g\left( y \right)}{x-y}}dy=\frac{1}{\pi }\int\limits_{E}{\frac{g{{\left( y \right)}^{2}}}{g\left( y \right)\left( x-y \right)}}dy \\ 
 & =\frac{1}{\pi }\int\limits_{E}{\left[ \sum\limits_{j=0}^{\ell +1}{{{\xi }_{j}}{{p}_{j}}\left( y \right)} \right]\frac{{{\rho }_{\ell -1}}\left( y \right)dy}{g\left( y \right)\left( x-y \right)}} \\ 
 & = \sum\limits_{j=0}^{\ell +1}\xi_j\int\limits_E  \left[\frac{1}{ x-y } p_j \left( y \right) \right]\omega_E\left( y \right)dy \\ 
 & =-\sum\limits_{j=0}^{\ell}{{{\xi }_{j+1}}{{q}_{j}}\left( x \right)}\quad\text{ for } x\in E.  
\end{aligned}\]
 On the other hand, we can write \[\tilde{g}={{\lambda }_{\ell }}{{p}_{\ell }}+{{\lambda }_{\ell -1}}{{p}_{\ell -1}}+\cdots +{{\lambda }_{1}}{{p}_{1}}+\lambda_0 \quad\text{ on } E\]
and it follows from the finite inversion of Hilbert transform that
\[\begin{aligned}
   \rho_{\ell -1}( x)
 & =\lambda_\ell  q_{\ell -1}\left( x \right)+\lambda_{\ell -1} q_{\ell -2}\left( x \right)+\cdots +\lambda_1 q_0\left( x \right)  \text{   for  } x\in E.
\end{aligned}\]
It follows from Theorem 4 that
\[\int\limits_{E}{{{\left| g\left( y \right) \right|}^{2}}{{\omega }_{E}}\left( y \right)}dy=\int\limits_{E}{{{\left| \tilde{g}\left( y \right) \right|}^{2}}{{\omega }_{E}}\left( y \right)}dy\]
so
\[\begin{aligned}
\sum\limits_{j=0}^{\ell }{{{\left| {{\lambda }_{j}} \right|}^{2}}}
&=-\int\limits_{E}{\left( x-{{a}_{1}} \right)\left( x-{{a}_{2}} \right)\cdots \left( x-{{a}_{2\ell }} \right){{\omega }_{E}}\left( x \right)dx}\\
&=\frac{1}{\pi }\int\limits_{E}{g\left( x \right){{\rho }_{\ell -1}}\left( x \right)dx.} \end{aligned}\]
Now we prove that if $\ell=1$ ($E$ is one compact interval) then
\[\int_E q_n(x)q_m(x)\frac{dx}{\pi^2\omega_E(x)}=\delta(n-m).
\]
Indded, it follows from the definition of $q_n$ that $-q_n=H(\pi p_{n+1}\omega_E)$ so apply the primitive  inversion
\[f(t)=\frac{1}{\pi \sqrt{(t-a)(b-t)}}\left( \int_{a}^{b}{\frac{Hf(s)}{s-t}}\sqrt{(s-a)(b-s)}ds+\int_{a}^{b}{f}(s)ds \right)\] for one interval 
 we have $p_{n+1}=H(q_n/\pi\omega_E)$. Apply theorem 4 we have
\[\int_E |q_n(x)|^2\frac{dx}{\pi^2\omega_E(x)}=\int_E |p_{n+1}(x)|^2\omega_E(x)dx=1.
\]
Moreover, it is easy to extend Theorem 4 to the form
\[\int_E f(x)\overline{g(x)}\omega_E(x)dx=\int_E \tilde f(x)\overline{\tilde g(x)}\omega_E(x)dx
\] and get the orthonormality of $q_n$ with respect to the positive measure $\frac{dx}{\pi^2\omega_E(x)}$ in the case $\ell=1$.
In general case we have
 \[p_n(x)\rho_{\ell-1}(x)=\frac 1\pi\int_E \frac{q_{n-1}(y)\rho_{\ell-1} (y)}{\pi\omega_E(y)}\frac{dy}{x-y} \quad\text{ for }\quad n\ge\ell.
\]
For example, let $E=\left[ -1,1 \right].$ Then 
\[g\left( x \right)=\left\{ \begin{matrix}
   \sqrt{1-{{x}^{2}}} & \text{if }\left| x \right|<1  \\
   0 & \text{otherwise,}  \\
\end{matrix} \right.\]
${{p}_{n}}=\sqrt{2}{{T}_{n}}$ and ${{\alpha }_{n}}=1/2$  for $n=1,2,\cdots .$ 
Here, \[{{T}_{n}}\left( x \right)=\cos n\theta =\frac{{{\left( x+\sqrt{{{x}^{2}}-1} \right)}^{n}}+{{\left( x-\sqrt{{{x}^{2}}-1} \right)}^{n}}}{2}\]
 is Chebisev polynomial of first kind. Moreover, ${{\alpha }_{0}}=1/\sqrt{2}$
 and ${{\beta}_{n}}=0$ for $n=0,1,2,\cdots .$ We have ${{q}_{0}}=\sqrt{2}=\sqrt{2}{{U}_{0}}$
  and 
\[\begin{aligned}
 {{q}_{1}}\left( y \right) & =\int\limits_{E}{\frac{{{p}_{2}}\left( x \right)-{{p}_{2}}\left( y \right)}{x-y}{{\omega }_{E}}\left( x \right)dx}\text{   } \\ 
 & =2\sqrt{2}\int\limits_{E}{\left( x+y \right){{\omega }_{E}}\left( x \right)dx} \\ 
 & =2\sqrt{2}y=\sqrt{2}{{U}_{1}}\left( y \right)  
\end{aligned}\]
where \[{{U}_{n-1}}\left( x \right)=\frac{\sin n\theta }{\sin \theta }=\frac{{{\left( x+\sqrt{{{x}^{2}}-1} \right)}^{n}}-{{\left( x-\sqrt{{{x}^{2}}-1} \right)}^{n}}}{2\sqrt{{{x}^{2}}-1}}\]  is Chebisev polynomial of second kind. 
But $q_n$'s satisfy the same linear recurrence as $U_n$'s so we have ${{q}_{n}}=\sqrt{2}{{U}_{n}}$  for every 
$n=0,1,2,\cdots .$ Thus, \[\frac{{{U}_{n-1}}\left( x \right)}{{{T}_{n}}\left( x \right)}\to {{\pi \tilde{\omega }}_{\left[ -1,1 \right]}}\left( x \right)\text{   as }n\to \infty \]
 for  $\left| x \right|>1.$  More exactly,
\[\begin{aligned}
 \frac{1}{\pi }\int\limits_{-1}^{1}{\frac{1}{\sqrt{1-{{y}^{2}}}}}\frac{dy}{x-y} & =\underset{n\to \infty }{\mathop{\lim }}\,\frac{{{U}_{n-1}}\left( x \right)}{{{T}_{n}}\left( x \right)} \quad\text{    if   }|x|>1\\ 
  & =\frac{1}{\sqrt{{{x}^{2}}-1}}\quad\text{    if   }x>1 \\ 
 & =-\frac{1}{\sqrt{{{x}^{2}}-1}}\quad \text{   if   }x<-1\\
  & =0\quad\text{    if   }|x|<1. \\
\end{aligned}\]
We also have \[\frac{1}{\pi }\int\limits_{-1}^{1}{\frac{{{T}_{n}}\left( x \right)-{{T}_{n}}\left( y \right)}{x-y}}\frac{dx}{\sqrt{1-{{x}^{2}}}}={{U}_{n-1}}\left( y \right)
\quad
\text{  for every } y\in \mathbb{C}.\]
On the other hand, $g{{\left( x \right)}^{2}}=1-{{x}^{2}}={{\xi }_{0}}+{{\xi }_{1}}{{p}_{1}}\left( x \right)+{{\xi }_{2}}{{p}_{2}}\left( x \right)$ so
 $\xi_1=0$ and ${{\xi }_{2}}=-\frac{1}{2\sqrt{2}}.$ Therefore, $\tilde{g}(x)=-{{\xi }_{2}}{{q}_{1}}\left( x \right)=x$ for $x\in E$.
Similarly, for the probability measure 
\[\frac{2\sqrt{1-{{x}^{2}}}}{\pi }dx\]  on $\left[ -1,1 \right]$ we have ${{p}_{n}}={{U}_{n}}$ and ${{q}_{n}}=2{{U}_{n}}.$ Thus,
\[ \begin{aligned}
 \frac{1}{\pi } \int_{-1}^1 \frac{\sqrt{1-{{y}^{2}}}}{x-y}dy&=\underset{n\to \infty }{\mathop{\lim }}\,\frac{{{U}_{n-1}}\left( x \right)}{{{U}_{n}}\left( x \right)}\quad\text{    if   }|x|>1 \\ 
 & =x-\sqrt{{{x}^{2}}-1}\quad\text{    if   }x>1 \\ 
 & =x+\sqrt{{{x}^{2}}-1}\quad\text{   if   }x<-1\\
   & =x\quad\text{    if   }|x|<1. \\ 
\end{aligned}\]
We also have  \[\frac{1}{\pi }\int\limits_{-1}^{1}{\frac{{{U}_{n}}\left( x \right)-{{U}_{n}}\left( y \right)}{x-y}}\sqrt{1-{{x}^{2}}}dx={{U}_{n-1}}\left( y \right) \quad \text{   for every }y\in \mathbb{C}.\]

\bigskip
\bigskip

\sect{\bf   	Logarithmic Potentials}

\bigskip
\par\noindent Let $w\left( x \right)={{e}^{-Q\left( x \right)}}>0$ satisfying $xw\left( x \right)\to 0$ as $x\to \pm \infty $. Consider the optimization problem
	\[\iint{\ln \frac{1}{\left| x-t \right|}}d\mu \left( x \right)d\mu \left( t \right)+2\int{Q\left( x \right)d\mu \left( x \right)}\to \min\]
subject to every Borel probability measure $\mu $ on the real line.  Let 
\[{{U}^{\mu }}\left( x \right)=\int{\ln \frac{1}{\left| x-t \right|}}d\mu \left( t \right)\]
 denote the potential of  $\mu $. Then
by \cite{Saff}
 there is a unique Borel probability measure ${{\mu }_{w}}$ on the real line solving this optimization problem and  ${{U}^{{{\mu }_{w}}}}\left( x \right)+Q\left( x \right)=F_w$ const for all $x\in\text{supp} \left( {{\mu }_{w}} \right)=:S_w$ and ${{U}^{{{\mu }_{w}}}}\left( x \right)+Q\left( x \right)\ge F_w$ for all $x\in \mathbb{R}.$ More exactly, ${{\mu }_{w}}$  is absolutely continuous and having compact support.  Let 
\[d{{\mu }_{w}}\left( t \right)=\frac{1}{\pi }f\left( t \right)dt\] 
then $F\left( x \right)+Q\left( x \right)$ is constant on $S_w $.
Here, 
\[F\left( x \right)=\frac{1}{\pi }\int\limits_{-\infty }^{\infty }{f\left( t \right)\ln \frac{1}{\left| x-t \right|}}\cdot dt=U^{\mu_w}(x).
\]  Moreover, 
\[\frac{1}{\pi }\int\limits_{S_w }{f\left( t \right)dt}=1.\] 
Hence,  ${F}'\left( x \right)+{Q}'\left( x \right)=0$ so 
$\tilde{f}\left( x \right)={Q}'\left( x \right)$  for $x\in S_w .$
It is proved in \cite{Saff} that if $Q$ is convex then $S_w$ is one interval. For example, if $w\left( x \right)={{W}_{\alpha }}\left( x \right)=\exp \left( -{{\left| x \right|}^{\alpha }} \right)$ is Freud weight then $S_w =\left[ -a\left( \alpha  \right),a\left( \alpha  \right) \right]$ if $\alpha >1$with
\[a\left( \alpha  \right)={{\left( \frac{2\alpha }{\pi }\int\limits_{0}^{\pi /2}{{{\sin }^{\alpha }}\vartheta d\vartheta } \right)}^{-1/\alpha }}={{\left( \frac{\sqrt{\pi }\Gamma \left( \alpha /2 \right)}{2\Gamma \left[ \left( \alpha +1 \right)/2 \right]} \right)}^{1/\alpha }}.\]
  If $\alpha =2m$ is an even integer then 
\[a\left( 2m \right)={{\left( \frac{4m}{\pi }\int\limits_{0}^{\pi /2}{{{\sin }^{2m}}\vartheta d\vartheta } \right)}^{-1/\left( 2m \right)}}=\sqrt[2m]{\frac{\left( 2m-2 \right)!!}{\left( 2m-1 \right)!!}}.\]
If  $Q$ is a polynomial then \[{{S}_{w}}=\bigcup\limits_{k=1}^{\ell }{\left[ {{a}_{2k-1}},{{a}_{2k}} \right]}\]  is a finite  union of intervals.  
Let 
\[K\left( x \right)=\prod\limits_{j=1}^{2\ell }{\left( x-{{a}_{j}} \right)}.
\]
 Then $K(x)\le 0$ for $x\in S_w $. The equilibrium measure of $S_w $ is 
\[d{{\nu }_{S_w}}\left( x \right)=\omega_{S_w}(x)dx= \frac{1}{\pi }\cdot \frac{\left| {{\rho }_{\ell -1}}\left( x \right) \right|}{\sqrt{\left| K\left( x \right) \right|}}\cdot dx\]
 where ${{\rho }_{\ell -1}}\left( x \right)={{x}^{\ell -1}}+\cdots $ is that unique polynomial satisfying 
\[\int\limits_{{{a}_{2j}}}^{{{a}_{2j+1}}}{\frac{{{\rho }_{\ell -1}}\left( x \right)}{\sqrt{K\left( x \right)}}}\cdot dx=0\]
 for $j=1,2,\cdots ,\ell -1.$ 
Clearly, $\omega_{S_w}\notin L^p $ for any $p>2$.  Moreover, 
\[\begin{aligned}
 {{F}_{w}} & =\int\limits_{{{S}_{w}}}{Q\left( x \right){{\omega }_{{{S}_{w}}}}\left( x \right)dx}+\int\limits_{{{S}_{w}}}{{{U}^{{{\mu }_{w}}}}\left( x \right){{\omega }_{{{S}_{w}}}}\left( x \right)dx} \\ 
 & =\int\limits_{{{S}_{w}}}{Q\left( x \right){{\omega }_{{{S}_{w}}}}\left( x \right)dx}+\frac{1}{\pi }\int\limits_{{{S}_{w}}}{f\left( t \right)dt}\int\limits_{{{S}_{w}}}{\ln \frac{1}{\left| t-x \right|}{{\omega }_{{{S}_{w}}}}\left( x \right)dx} \\ 
 & =\int\limits_{{{S}_{w}}}{Q\left( x \right){{\omega }_{{{S}_{w}}}}\left( x \right)dx}-\ln \text{cap}\left( {{S}_{w}} \right).  
\end{aligned}\]
If $\ell=1$ we have the density 
\[f\left( x \right)=\frac{1}{\sqrt{\left| K\left( x \right) \right|}}\left( 1+\int\limits_{{{S}_{w}}}{\frac{{Q}'\left( t \right)\sqrt{\left| K\left( t \right) \right|}}{t-x}\cdot dt} \right).\]
If we know the density $f\in {{L}^{p}}$ with $p>2$ then 
\[f\left( x \right)=\frac{1}{\pi {{\omega }_{S_w}}\left( x \right)}\int\limits_{{{S}_{w}}}{\frac{{Q}'\left( y \right){{\omega }_{S_w}}\left( y \right)}{y-x}}dy.\]
For example, if $w\left( x \right)={{W}_{\alpha }}\left( x \right)=\exp \left( -{{\left| x \right|}^{\alpha }} \right)$ then 
\[\begin{aligned}
 f\left( x \right) & =\frac{1}{\sqrt{{{a}^{2}}-{{x}^{2}}}}\left( 1+2\alpha \int\limits_{0}^{a}{\frac{{{t}^{\alpha }}\sqrt{{{a}^{2}}-{{t}^{2}}}}{{{t}^{2}}-{{x}^{2}}}dt} \right) \\ 
 & =\frac{2\alpha \sqrt{{{a}^{2}}-{{x}^{2}}}}{\pi }\int\limits_{0}^{a}{\frac{{{t}^{\alpha }}dt}{\left( {{t}^{2}}-{{x}^{2}} \right)\sqrt{{{a}^{2}}-{{t}^{2}}}}}.  
\end{aligned}\]
The last identity holds in the case  $f\in {{L}^{p}}$ for some $p>2$ only. If $w\left( x \right)=\exp \left( -{{x}^{2}} \right)$ is the Gaussian then ${{S}_{w}}=\left[ -1,1 \right]$ and 
\[d{{\mu }_{w}}\left( x \right)=\frac{2\sqrt{1-{{x}^{2}}}}{\pi }dx\]
 is the optimizing measure. More generally, if $w\left( x \right)=\exp \left( -{{x}^{2m}} \right)$  then ${{S}_{w}}=\left[ -a,a \right]$ with
\[a=\sqrt[2m]{\frac{\left( 2m-2 \right)!!}{\left( 2m-1 \right)!!}}\]
and
\[d{{\mu }_{w}}\left( x \right)=\frac{2m\sqrt{{{a}^{2}}-{{x}^{2}}}}{\pi }\sum\limits_{k=0}^{m-1}{\frac{\left( 2k -1 \right)!!}{\left( 2k \right)!!}}{{a}^{2k}}{{x}^{2m-2k -2}}dx\]
 is the optimizing measure. Here $(-1)!!=0!!=1$. Now we focus our attention on the conductor $\Sigma =\left[ -1,1 \right]$ and consider the optimization problem
	\[\iint{\ln \frac{1}{\left| x-t \right|}}d\mu \left( x \right)d\mu \left( t \right)+2\int{Q\left( x \right)d\mu \left( x \right)}\to \min\]
subject to every Borel probability measure $\mu$ supported in the conductor $\Sigma =\left[ -1,1 \right]$. There is exactly one measure ${{\mu }_{w}}$ solves this optimization problem.  If $Q=\varepsilon {{T}_{n}}$ (${{T}_{n}}$ denotes the Chebisev polynomial of the first kind) then $\text{supp}\left( {{\mu }_{w}} \right)=\left[ -1,1 \right]$ and 
\[d{{\mu }_{w}}\left( x \right)=\frac{1-n\varepsilon {{T}_{n}}\left( x \right)}{\pi \sqrt{1-{{x}^{2}}}}dx\] provided that $\left| \varepsilon  \right|\le 1/n.$
In fact, the Hilbert transform $\tilde{f}$ of the density function $f$ should satisfy $\tilde{f}={Q}'=n\varepsilon {{U}_{n-1}}$  (${{U}_{n-1}}$ denotes the Chebisev polynomial of the second kind) so by the inversion of finite Hilbert transform we have 
\[f\left( x \right)=\frac{1-n\varepsilon {{T}_{n}}\left( x \right)}{\sqrt{1-{{x}^{2}}}}\] 
which is certainly positive if $\left| \varepsilon  \right|\le 1/n.$ More generally, let $Q=\sum\limits_{k=0}^{n}{{{\varepsilon }_{k}}{{T}_{k}}}$ with $\sum\limits_{k=1}^{n}{k\left| {{\varepsilon }_{k}} \right|\le 1.}$ Then $\text{supp}\left( {{\mu }_{w}} \right)=\left[ -1,1 \right]$ and
\[d{{\mu }_{w}}\left( x \right)=\frac{1-\sum\limits_{k=1}^{n}{k{{\varepsilon }_{k}}{{T}_{k}}\left( x \right)}}{\pi \sqrt{1-{{x}^{2}}}}dx.\]
Now let $Q(x)=-\epsilon T_2(x)$ with $\epsilon>1/2$. Let $a=4\epsilon>2$. We have $\tilde f(x)=Q'(x)=-2ax$ and $S_w=[-1,-\alpha]\cup[\alpha,1]$ with $\alpha=\sqrt{\frac{a-2}a}$. Moreover, \[d\mu_w(x)=\frac {\sqrt{a}|x|}{\pi}\sqrt{\frac{ax^2+2-a}{1-x^2}}dx
\] 
Similarly, if $Q(x)=\epsilon T_2(x)$ with $\epsilon>1/2$ then $S_w=[-1/\sqrt a, 1/\sqrt a]$ and $d\mu_w(x)=2a\sqrt{{1-ax^2}}dx/\pi$ ($a=4\epsilon>2$).
If $Q(x)=\epsilon T_4(x)$ with $\epsilon>1/4$ then $S_w=[-\beta,-\alpha]\cup[\alpha,\beta]$ with $$\alpha^2=\frac12-\frac1{4\sqrt\epsilon},\quad\beta^2=\frac12+\frac1{4\sqrt\epsilon}\quad\text{and}\quad  
 \frac{d\mu_w(x)}{dx}=\frac{2|x|[1-4\epsilon(2x^2-1)^2]}{\pi\sqrt{(x^2-\alpha^2)(\beta^2-x^2)}} .$$
If $Q(x)=-\epsilon T_4(x)$ with $\epsilon>1/4$ then $S_w=[-1,-\alpha]\cup[-\beta,\beta]\cup[\alpha,1]$ is of 3 intervals.
If  $Q=\sum\limits_{k=0}^{n}{{{\varepsilon }_{k}}{{T}_{k}}}$  is a polynomial of degree $n$ then the number of intervals of  $S_w $ is at most $1+ n$. It can be proved easily by iterated balayage algorithm (next section). Moreover, if ${1-\sum\limits_{k=1}^{n}{k{{\varepsilon }_{k}}{{T}_{k}}\left( x \right)}}$ is increasing on $[-1,1]$ then $S_w$ is one interval containing 1 and if ${1-\sum\limits_{k=1}^{n}{k{{\varepsilon }_{k}}{{T}_{k}}\left( x \right)}}$ is decreasing on $[-1,1]$ then $S_w$ is one interval containing $-1$. If 
${1-\sum\limits_{k=1}^{n}{k{{\varepsilon }_{k}}{{T}_{k}}\left( x \right)}}$ is decreasing on $[-1,t_0]$ 
and increasing on $[t_0,1]$
then $S_w$ is one interval containing $-1$ or one interval containing $1$ or a union of one interval containing $-1$ with one interval containing $1$.
\bigskip
\bigskip

\sect{\bf   	Balayage onto a compact set }

\bigskip
\par\noindent Let $\mathfrak{K}$ be  a compact subset of the complex plane of positive logarithmic capacity such that $\bar{\mathbb{C}}\backslash \mathfrak{K}$ is regular for Dirichlet problem. Let $\nu $ be a positive Borel measure of compact support on the complex plane. Then there is a unique positive measure $\hat{\nu }$ supported on $K$ such that $\left\| {\hat{\nu }} \right\|=\left\| \nu  \right\|$ and ${{U}^{{\hat{\nu }}}}-{{U}^{\nu }}$ is constant on  $\mathfrak{K}$. The measure $\hat{\nu }$ is called the balayage of $\nu $ onto $\mathfrak{K}$ and denoted by $\text{Bal}\left( \nu ,\mathfrak{K} \right)$. For a signed measure $\sigma ={{\sigma }^{+}}-{{\sigma }^{-}}$ we define $\text{Bal}\left( \sigma ,\mathfrak{K} \right)=\text{Bal}\left( {{\sigma }^{+}},\mathfrak{K} \right)-\text{Bal}\left( {{\sigma }^{-}},\mathfrak{K} \right).$ For example, let  \[E=\bigcup\limits_{k=1}^{\ell }{\left[ {{a}_{2k-1}},{{a}_{2k}} \right]}\]   and    \[K\left( x \right)=\prod\limits_{j=1}^{2\ell }{\left( x-{{a}_{j}} \right)}.\] Fix a point $s$ in a gap of $E$. Then the density function of the  balayage of the point mass ${{\delta }_{s}}$ onto $E$ is
\[\frac{d}{dt}\text{Bal}\left( {{\delta }_{s}},E \right)=\frac{1}{\pi }\sqrt{\frac{K\left( s \right)}{\left| K\left( t \right) \right|}}\cdot \left| \frac{R _{s}^{*}\left( t-s \right)}{t-s} \right|\quad \text{  for }t\in E.\]
Here, $R_s\left( t \right)={{t}^{\ell -1}}+\cdots$ is a monic polynomial of degree $\ell-1$ 
without repeated root which is 
uniquely determined by $s$ and $E$   and  $R_s^{*}\left( t \right)=t^{\ell -1}R_s\left( 1/t \right)$ is the reciprocal polynomial of $R_s.$  Moreover, in each gap of $(E-s)^{-1}$
there is exactly one root of $R_s$. 
In fact, the balayage of the point mass ${{\delta }_{s}}$ onto $E$ is exactly the image of the equilibrium measure of \[(E-s)^{-1}:=\{(t-s)^{-1}:\text{ } t\in E\}\] 
under the mapping $(t-s)^{-1}\to t$ which maps  (one-to-one) $(E-s)^{-1}$ onto $E$.
Let $P\left( t \right)=R_{s}^{*}\left( t-s \right)$.  Then $P$ is of degree $ \ell -1$ or $\ell -2$. Moreover, in each gap of $E$ not including $s$ there is exactly one root of $P.$  If $P$ is of degree $\ell -1$ then there is one more zero in $\mathbb{R}\backslash \left[ {{a}_{1}},{{a}_{2\ell }} \right].$
Let $\nu $ be a finite positive Borel measure supported in gaps of $E$. Then the density function of the balayage of $\nu $ onto $E$ is 
\[\frac{d}{dt}\text{Bal}\left( \nu ,E \right)=\frac{1}{\pi \sqrt{\left| K\left( t \right) \right|}}\sum\limits_{k=1}^{\ell -1}{\int\limits_{{{a}_{2k}}}^{{{a}_{2k+1}}}{\left| \frac{R_{s}^{*}\left( t-s \right)}{t-s} \right|\sqrt{K\left( s \right)}d\nu \left( s \right),}}\quad t\in E.\]
Here we use methods of \cite{Damelin} and \cite{Saff} to get these explicit formulas for balayages.
For example, let  $E$ be a compact subset of the real line and $\nu$  a Borel positive finite measure of compact support such that  $\nu(E)=0$.   If   $E=[a,b]$ then Bal$(\nu,E)$   is absolutely continuous and $$\frac{d}{dt}\text{Bal}(\nu,E)=\frac 1{\pi\sqrt{(t-a)(b-t)}}\int\frac{\sqrt{(s-a)(b-s)}}{|s-t|}d\nu(s).$$ 
   Thus, if supp$(\nu)\subseteq(-\infty,a]$  then ${\pi\sqrt{(t-a)(b-t)}}\frac d{dt}\text{Bal}(\nu,E)$  is increasing in  $[a,b]$. Similarly, if  supp$(\nu)\subseteq[b,\infty)$  then ${\pi\sqrt{(t-a)(b-t)}}\frac d{dt}\text{Bal}(\nu,E)$  is decreasing in  $[a,b]$.
If   $E=[a,b]\cup[c,d]$ and supp$(\nu)\subseteq[b,c]$  then  $${\pi\sqrt{(t-a)(d-t)}}\frac d{dt}\text{Bal}(\nu,E)$$  is increasing in$[a,b]$ and decreasing  in $[c,d]$.
Now we return to the external field $Q=\sum\limits_{k=0}^{n}\varepsilon_kT_k$   a polynomial of degree $n$ with $\sum\limits_{k=0}^{n}k|\varepsilon_k|>1$ on the conductor $\Sigma=[-1,1]$. We let $\ell(E)$ denote the number of intervals of $E$ if $E$ is a finite union of intervals. Let 
 \[f_0(x)=\frac{1-\sum\limits_{k=1}^{n}{k{{\varepsilon }_{k}}{{T}_{k}}\left( x \right)}}{\sqrt{1-{{x}^{2}}}}\]  and
\[E_1=\bigl\{x\in[-1,1]:\text{ }f_0(x)\ge0 \bigr\}.\]
Then $\ell(E_1)\le[n/2]+1<n+1$. Let $f_1(x)$ denote the density function of the balayage of  the signed measure$f_0(x)dx$ onto $E_1$. Then $\tilde f_1=Q'$ so
\[f_1\left( x \right)=\frac{1}{\pi g_{E_1}\left( x \right) } \left( \int\limits_{E_1}\frac{\tilde{g}_{E_1}\left( x \right)-\tilde{g}_{E_1}\left( y \right)}{x-y} f_1\left( y \right)dy+\int\limits_{{E_1}}\frac{g_{E_1}\left( y \right)Q'\left( y \right)}{y-x}dy \right)\]
Thus, $f_1g_{E_1}$ is a polynomial of degree $\le n+\ell(E_1)-1$ so the set $E_2=\bigl\{x\in E_1:\text{ }f_1g_{E_1}(x)\ge0 \bigr\}$ is a union of  $\ell(E_2)\le n+1$ intervals. Let $f_2(x)$ denote the density function of the balayage of  the signed measure $f_1(x)dx$ onto $E_2$. Then $\tilde f_2=Q'$ so $f_2g_{E_2}$ is a polynomial of degree $\le n+\ell(E_2)-1$. Repeat this argument infinite times we get a sequence $f_0,f_1,\cdots$ converging to a positive function $f=P/g_E$ where $E$ is a union of $\ell\le n+1$ intervals and $P$ is polynomial of degree $\le n+\ell-1$. Clearly, $S_w=E$ and $d\mu_w(x)=f(x)dx/\pi$. On the other hand, if
$\sqrt{1-x^2}f_0(x)={1-\sum\limits_{k=1}^{n}{k{{\varepsilon }_{k}}{{T}_{k}}\left( x \right)}}$ is increasing in $[-1,1]$ then $E_1$ is one interval containing 1. Therefore, every $E_m$ is an interval containing 1 so $S_w$ itself is one interval containing 1. Similarly, if $\sqrt{1-x^2}f_0(x)={1-\sum\limits_{k=1}^{n}{k{{\varepsilon }_{k}}{{T}_{k}}\left( x \right)}}$ is decreasing in $[-1,1]$ then $E_1$ is one interval containing $-1$. Therefore, every $E_m$ is an interval containing $-1$ so $S_w$ itself is one interval containing $-1$.
Now we are interested in the inversion problem of balayages. For example, let $E=[0,1]$. We  look after the probability measure $\nu$ supported in  $[2,3]$ such that the balayage of $\nu$ onto $E$ is the equilibrium measure of $E$. Let  $\varphi$ denote the density function of $\nu$ then
\[\int_2^3 \frac { \sqrt{s(s-1)}\varphi(s)}{s-t}ds=1\qquad\text{ for } t\in[0,1].
\]
Moreover,
\[\int_2^3 \varphi(s)\ln\frac 1{|t-s|}ds=\text{const}\qquad\text{ for } t\in[0,1].
\]
Taking derivative two times according to $t$ we have
\[\int_2^3 \frac {\varphi(s) }{|t-s|^2}ds=0\qquad\text{ for } t\in[0,1],
\]
which is impossible bacause  $\varphi$ is positive. This means the equilibrium measure of a compact set cannot be  the balayage of a probabilty measure supported outside of the compact set.
\bigskip
\bigskip

\sect{\bf   Singular Integral Equations }

\bigskip
\par\noindent 
Several authors \cite{Cooke}  \cite{ESTRADA}  \cite{Garoufalidis} \cite{Gautesen} 
 \cite{MARGETSON}  \cite{Morland} 
 study the logarithmic integral equation (for water waves, random matrices, etc. )
\[F\left( x \right)=\frac{1}{\pi }\int\limits_{E}{f\left( y \right)\ln \frac{1}{\left| x-y \right|}dy},
\quad x\in E, 
\]
where $E$ is a finite union of compact intervals and $F$ is smooth on $E.$ It follows at one from the condition of equilibrium measure that
\[\int\limits_{E}{f\left( y \right)dy}=-\frac{\pi }{\ln \text{cap}\left( E \right)}\int\limits_{E}{F\left( y \right){{\omega }_{E}}\left( y \right)dy} \quad\text{ if cap}(E)\ne 1.
\]
If $\text{cap}\left( E \right)=1$ then 
\[\int\limits_{E}{F\left( y \right){{\omega }_{E}}\left( y \right)dy}=0.\]
Theorem 3 shows that if $F$ is absolutely continuous with  ${F}'\in {{L}^{p}}$ $\left( p>2 \right)$ then this integral equation has at most one solution $f \in {{L}^{p}}$  determined by the explicit formula
\[f\left( x \right)=\frac{1}{\pi {{\omega }_{E}}\left( x \right)}\int\limits_{E}{\frac{{F}'\left( y \right){{\omega }_{E}}\left( y \right)}{x-y}}dy\quad \text{     for a}\text{.e}\text{. }x\in E.\]
In fact, the weak derivative of $-F$ is exactly the finite Hilbert transform of $f$.
Here, we do not need the smoothness of function $F$ as authors have requested to solve this equation. Moreover, if $F$ is non-zero constant (infinitely differentiable) then this equation has no solution in 
$L^p$ for any $p>2$.  If $E=[a,b]$  ($b-a\ne 4$)
is a compact interval we need only the weak derivative of $F$ belonging to $L^p$ with $p>1$ and the solution in $L^p$   is determined uniquely  by formula  \cite{Chakrabarti1}  \cite{Cooke}
\[f(t)=\frac{1}{ \sqrt{(t-a)(b-t)}}\left[\frac 1\pi \int_{a}^{b}{\frac{{F}'(s)}{t-s}}\sqrt{(s-a)(b-s)}ds\right.\]
\[\left.+
\left(\ln\frac 4{b-a}\right)^{-1}
\int_{a}^{b}\frac{F(s)ds}{\sqrt{(s-a)(b-s)}} \right].\]
If $E=\left[ -b,-a \right]\cup \left[ a,b \right]$ with ${{b}^{2}}-{{a}^{2}}\ne 4$ and $F$ is even then
\[f\left( x \right)=\frac{2|x|}{ \sqrt{\left( {{b}^{2}}-{{x}^{2}} \right)\left( {{x}^{2}}-{{a}^{2}} \right)}}\left[\left(\ln \frac{4}{b^2-a^2}\right)
^{-1}
\int\limits_{a}^{b}{\frac{F\left( y \right)ydy}{\sqrt{\left( {{b}^{2}}-{{y}^{2}} \right)\left( {{y}^{2}}-{{a}^{2}} \right)}}}
\right.\]
\[\left.+\frac 1\pi \int\limits_{a}^{b}{\frac{\sqrt{\left( {{b}^{2}}-{{y}^{2}} \right)\left( {{y}^{2}}-{{a}^{2}} \right)}{F}'\left( y \right)dy}{{{x}^{2}}-{{y}^{2}}}} \right].\]
Manam \cite{Manam} 
 studied the logarithmic integral equation
\[\frac{1}{\pi }\int\limits_{E}{f\left( y \right)\ln \left| \frac{x+y}{x-y} \right|dy}=G\left( x \right),
\quad x\in E, 
\]
where $E$ is a finite union of positive compact intervals and $G$ is smooth on $E.$ 
If ${G}'\in {{L}^{p}}$ ($p>2$), Theorem 3 shows that there is at most one solution in $L^p$ determined by the explicit formula 
\[f\left( x \right)=\frac{2}{\pi {{\omega }_{{{E}^{2}}}}\left( {{x}^{2}} \right)}\int\limits_{E}{\frac{{G}'\left( y \right){{\omega }_{{{E}^{2}}}}\left( {{y}^{2}} \right)ydy}{{{x}^{2}}-{{y}^{2}}}}\quad \text{     for a}\text{.e}\text{. }x\in E.\]
Here, $E^2=\{x^2:\text{ } x\in E\}$.
In fact, $\tilde{f}\left( t \right)+\tilde{f}\left( -t \right)=-{G}'\left( t \right)$ so we have
\[\frac{1}{\pi }\int\limits_{{{E}^{2}}}{\frac{f\left( \sqrt{y} \right)dy}{y-t}}={G}'\left( \sqrt{t} \right),
\qquad t\in E^2
\]
and Theorem 3 is applied to get the function $f$ uniquely in $L^p$ ($p>2$). 
If $G$ is identically non-zero constant then this logarithmic equation has no solution in $L^p$ for any $p>2$. 
 Moreover, if $E=[a,b]$
is a compact interval then
 \[f(t)=\frac{2}{\pi \sqrt{({{t}^{2}}-{{a}^{2}})({{b}^{2}}-{{t}^{2}})}}\left( \int_{a}^{b}{\frac{s{G}'(s)}{{{t}^{2}}-{{s}^{2}}}}\sqrt{({{s}^{2}}-{{a}^{2}})({{b}^{2}}-{{s}^{2}})}ds+\int_{a}^{b}{sf}(s)ds \right)\]
provided that ${G}'\in {{L}^{p}}$ with $p>1.$ 
Specially, if $E=[0,a]$  then  we have the unique solution \cite{Cooke}
\[f\left( x \right)=-\frac{2}{\pi }\frac{d}{dx}\int\limits_{x}^{a}{\frac{\alpha S\left( \alpha  \right)d\alpha }{\sqrt{{{\alpha }^{2}}-{{x}^{2}}}}},\]
where
\[S\left( \alpha  \right)=\frac{1}{\alpha }\frac{d}{d\alpha }\int\limits_{0}^{\alpha }{\frac{xG\left( x \right)dx}{\sqrt{{{\alpha }^{2}}-{{x}^{2}}}}}
=\frac{1}{\alpha }\frac{d}{d\alpha }\int\limits_{0}^{\alpha }\sqrt{\alpha ^2-x^2}G'(x)dx
.\]     
(Note that formulas (25) and (28) in \cite{Cooke} are incorrect.)
In fact, using the formula
\[\frac{1}{2}\ln \left| \frac{x+y}{x-y} \right|=\int\limits_{0}^{\min \left( x,y \right)}{\frac{tdt}{\sqrt{\left( {{t}^{2}}-{{x}^{2}} \right)\left( {{t}^{2}}-{{y}^{2}} \right)}}}\]
we have
\[G\left( x \right)=\frac{2}{\pi }\int\limits_{0}^{x}{\frac{S\left( t \right)tdt}{\sqrt{{{x}^{2}}-{{t}^{2}}}}}\quad
\text{    with  }S\left( \alpha  \right)=\int\limits_{\alpha }^{a}{\frac{f\left( t \right)dt}{\sqrt{{{t}^{2}}-{{\alpha }^{2}}}}}\]  
and apply inversion formulas of Abel integrals we have the unique solution.
Now we consider the following integral equation \cite{Bruckner}
\[\frac{1}{\pi }\int\limits_{0}^1{f\left( t \right)\ln \frac{1}{\left| x-t \right|}dt}=F\left( x \right),
\quad x\in  [2,3].
\]
Assume that the weak derivative ${F}'\in {{L}^{p}}$  $\left( p>1 \right)$. Let 
\[\phi(x)=\frac{1}{\pi \sqrt{(x-2)(3-x)}}\left( \int_2^3\frac{F'(s)}{s-x}\sqrt{(s-2)(3-s)}ds+
\int_0^1 f(t)dt\right)
\]
for  $ x\in  [2,3].$ Then
 \[\begin{aligned}
\int\limits_2^3\phi\left( y \right)dy
&=\int_0^1 f(t)dt\\
F\left( x \right)
&=\frac 1\pi \int_2^3 \phi\left( y \right) \ln \frac{1}{| x-y|}dy+\text{const,}\quad x\in  [2,3]\\
\end{aligned}
\]
so $\phi(x)dx$ is the balayage of $f(t)dt$ into $[2,3]$ and we have
\[\phi(x)=\frac{1}{\pi \sqrt{(x-2)(3-x)}}\int_{0}^{1}\frac{f(t)\sqrt{(t-2)(t-3)}}{x-t}dt \quad
\text{ for } x\in  [2,3].\]
This is also a singular integral equation which is not studied enough in literature. Only numerical simulations are made for approximate solution. 

\bigskip
\bigskip

\sect{Hilbert transform on positive semi-axis and water waves}

\bigskip
\par\noindent The singular integral equations in theory of water waves \cite{Banerjea} 
\cite{Chakrabarti}
request us to study the inversion of Hilbert transform 
\[\tilde{f}\left( t \right)=\frac{1}{\pi }\int\limits_{0}^{\infty }{\frac{f\left( \xi  \right)}{t-\xi }}d\xi\quad \text{    for   }t>0.\]
Let $\phi \left( t \right)=f\left( {{t}^{2}} \right)\text{sign}\left( t \right)$ be an odd function on the real line. Assume that $\phi\in L^p$ for some $p>1$. Then 
\[\int_0^\infty \frac{|f(\xi)|^pd\xi}{\sqrt\xi}<\infty\]
and the Hilbert transform of $\phi$ is an even function determined by the explicit formula
$\tilde\phi(x)=\tilde f(x^2)$ and we get the inversion formula
\[f\left( t \right)=\frac{\sqrt t}{\pi }\int\limits_{0}^{\infty }\frac{\tilde f\left( \xi  \right)d\xi}{(\xi -t)\sqrt\xi}\quad \text{    for   }t>0.\]
Now we consider the following singular integral equation appeared  in theory of water waves \cite{Banerjea} \cite{Chakrabarti}
\[\frac 1\pi \int_0^\infty f(t)\left[c\ln\frac{|x-t |}{|x+t|}+\frac1{x+t}+\frac1{x-t}
\right]dt=G(x)\]
 for  $x\in E:=[0,a]\cup [b,\infty)$  and $f$ is supported in $E$. The function $G$ is also known in $E$ only. 
Let
\[\lambda (t)=c\int_0^t f(\xi)d\xi +f(t)  \quad\text{ for }  t>0.\]
We get at once
\[\frac1\pi\int_0^\infty \frac{2{x}\lambda (t)dt}{x^2-t^2} = {G(x)} \quad\text{ for }  x>0.\]
Let $\phi \left( x \right)=\lambda \left( \left| x \right| \right)$. Then
\[\tilde{\phi }\left( x \right)=\frac{1}{\pi }\int\limits_0^\infty \lambda \left( t \right)\left( \frac{1}{x-t}+\frac{1}{x+t} \right)dt =\frac{1}{\pi }\int\limits_0^\infty \frac{2x\lambda \left( t \right)}{x^2-t^2}dt =G\left( x \right)\]
for $x>0$ and $\tilde{\phi }\left( x \right)=-G(-x)$ for $x<0$. We need only that $G\in L^p$ for some  $p>1$. Thus,
\[\begin{aligned}
  \lambda \left( x \right)& =-\frac{1}{{{\pi }^{}}}\int\limits_{0}^{\infty }{G\left( t \right)\left( \frac{1}{x-t}-\frac{1}{x+t} \right)dt} \\ 
 & =\frac{1}{{{\pi }^{}}}\int\limits_{0}^{\infty }{\frac{2tG\left( t \right)}{{{t}^{2}}-{{x}^{2}}}dt}.  
\end{aligned}\]
and
\[\begin {aligned}
f\left( x \right)
&=\frac{d}{dx}\left[ {{e}^{-cx}}\int\limits_{0}^{x}{{{e}^{c\xi }}}\left( \frac{1}{{{\pi }^{}}}\int\limits_{0}^{\infty }{\frac{2tG\left( t \right)}{{{t }^{2}}-{{\xi}^{2}}}dt} \right)d\xi  \right]\\
&=\lambda(x)-ce^{-cx}\int_0^x\lambda(\xi)e^{c\xi}d\xi.
\end {aligned}\]
But $f\left( x \right)=0$ for $x\in \left( a,b \right)$ so 
\[\lambda(x)=c\int_0^a f(\xi)d\xi \qquad\text{ for } x\in \left( a,b \right)\]
and consequently,
\[\frac{1}{\pi }\int\limits_{{{a}^{2}}}^{{{b}^{2}}}{\frac{G\left( \sqrt{t} \right)}{x-t}dt}=-c\int_0^af(\xi)d\xi- \frac{1}{\pi }\int\limits_E \frac{2\xi G\left( \xi  \right)}{x-\xi ^2}d\xi ,\quad\text{ for } x\in [a^2,b^2].\]
Therefore, for $t\in [a,b]$
\[G\left( t \right)=\frac{2}{\pi \sqrt{\left( t^2-a^2 \right)\left( b^2-t^2 \right)}}\left\{ \int\limits_a^b\xi G\left( \xi  \right)d\xi \right.\]
\[\left.-\int\limits_a^b \left[ c\int_0^a f(\xi)d\xi+\frac{1}{\pi }\int\limits_E\frac{2\xi G\left( \xi  \right)}{x^2-\xi ^2 } d\xi  \right]\frac{\sqrt{\left( x^2-a^2 \right)\left( b^2-x^2 \right)}}{x^2-t^2}xdx \right\}.\]
Here, the integrals $\int\limits_a^b\xi G\left( \xi  \right)d\xi $ and $\int_0^a f(\xi)d\xi $ are arbitrary constants.

\bigskip
\bigskip
\sect{AEROFOIL THEORY}

\bigskip
\noindent
Porter \cite{Porter} studied the following integro-differential equation appeared in aerofoil theory
\[\frac{1}{\pi }\sqrt{\frac{x-1}{x}}\int\limits_{1}^{\infty }{\sqrt{\frac{t}{t-1}}}\frac{{f}'\left( t \right)dt}{t-x}=\lambda f\left( x \right)+2\alpha \left( 1-\sqrt{\frac{x-1}{x}} \right)\text{ for   } x>1,\]
where $-2f\left( x \right)/\lambda $ is the slope of the jet and $\alpha ,\lambda $ are known parameters. Let
$\varphi \left( x \right)=f\left( 1/x \right)$. Then 
\[{f}'\left( x \right)=-\frac{{\varphi }'\left( 1/x \right)}{{{x}^{2}}}\] and 
\[-\frac{x\sqrt{1-x}}{\pi }\int\limits_{0}^{1}{\frac{t}{\sqrt{1-t}}}\frac{{\varphi }'\left( t \right)dt}{x-t}=\lambda \varphi \left( x \right)+2\alpha \left( 1-\sqrt{1-x} \right)\text{  for  } 0<x<1\]
or equivalently, 
\[\frac{1}{\pi }\int\limits_{0}^{1}{\frac{t{\varphi }'\left( t \right)}{\sqrt{1-t}}}\frac{dt}{x-t}=-\frac{\lambda \varphi \left( x \right)+2\alpha \left( 1-\sqrt{1-x} \right)}{x\sqrt{1-x}}\text{,    }\quad 0<x<1.\]
If we can write
\[\varphi \left( t \right)=\int\limits_{a}^{b}{u\left( ts \right)v\left( s \right)ds}=\int\limits_{a}^{b}{u\left( t\xi  \right)v\left( \xi  \right)d\xi }\]
with ${u}'\left( t \right)=1/\sqrt{t^3}$ and $v$ is continuous then
\[t{\varphi }'\left( t \right)=t\int\limits_{a}^{b}{{u}'\left( ts \right)v\left( s \right)sds}=\int\limits_{a}^{b}{\frac{v\left( s \right)ds}{\sqrt{ts}}}\]
and consequently, 
\[\frac{1}{\pi }\int\limits_{0}^{1}{\frac{t{\varphi }'\left( t \right)}{\sqrt{1-t}}}\frac{dt}{x-t}=\frac{1}{\pi }\int\limits_{a}^{b}{\frac{v\left( s \right)ds}{\sqrt{s}}}\int\limits_{0}^{1}{\frac{1}{\sqrt{\left( 1-t \right)t}}}\cdot \frac{dt}{x-t}=0\]
for $0<x<1.$  Therefore, 
\[\varphi \left( x \right)=-\frac{2\alpha \left( 1-\sqrt{1-x} \right)}{\lambda }\] is the unique solution of the form 
\[\varphi \left( t \right)=\int\limits_{a}^{b}{u\left( ts \right)v\left( s \right)ds}.\] 
The slope of jet  is \[\frac{4\alpha }{{{\lambda }^{2}}}\left( 1-\sqrt{\frac{x-1}{x}} \right).\]

\bigskip
\bigskip

\sect{Analytic Matrix Models and their Planar Limits}

\bigskip
\noindent
An admissible potential is a lower-semicontinuous function $V:\R\to\R$ satisfying \[
\lim_{|x|\to\infty}\frac{V(x)}{2\ln|x|}>1.\]
For an analytic random matrix model  \cite{Garoufalidis} with admissible potential $V$ 
we defined the planar limit 
\[
I^V=\inf_\mu I^V(\mu)=\inf_\mu \iint \ln\frac 1{|x-y|} d\mu(x)d\mu(y)+\int V(x)d\mu(x),
\]
where $\mu$ is running in the set of probability measures supported in $\R.$ It is well known that there is a unique probabilty measure $\mu^V$ such that $I^V=I^V(\mu^V)$.
A $1-$cut potential is an admissible potential $V$ such that the support of $\mu^V$ is a single interval $[-2c+b,2c+b]$. Then for $1-$cut potential $V$ we have the planar limit (the $F_w$ constant with $Q=V/2$)
\[\begin{aligned}
I^V&=\frac 1{2\pi}\int_{-2c+b}^{2c+b}\frac{V(x)dx}{\sqrt{4c^2-(x-b)^2}}-\ln c\\
&=\frac 1{\pi}\int_{-1}^{1}\frac{V(2cx+b)dx}{\sqrt{1-x^2}}-\ln c.\\
\end{aligned}
\]
Here, we do not need the smoothness of $V$ as  formula (28) of \cite{Garoufalidis} had requested.
The density function $f$ of $\mu^V$ will satisfying\[
\int_{-2c+b}^{2c+b}f(y)\ln\frac 1{|x-y|}dy=I^V-\frac{V(x)}2
\]
so\[
f(x)=\frac 1{\pi\sqrt{4c^2-(x-b)^2}}\left[1+\int_{-2c+b}^{2c+b}\frac{V'(y)\sqrt{4c^2-(y-b)^2}dy}{2(y-x)}\right]
\]
or equivalently,
\[
f(2cx+b)=\frac 1{2c\pi\sqrt{1-x^2}}\left[1+c\int_{-1}^{1}\frac{V'(2cy+b)\sqrt{1-y^2}dy}{y-x}\right].
\]
Here, we do not need the smoothness of $V$, only the local absolute continuity of $V$ with the weak derivative $V'\in L^p_{loc}$ for some $p>1$. Moreover, $c>0$ and $b\in\R$ maximize the function
\[\ln c- \frac 1{2\pi}\int_{-1}^{1}\frac{V(2cx+b)dx}{\sqrt{1-x^2}}.
\]
Taking derivative according to $c$ we have
\[\frac c\pi\int_{-1}^{1}\frac{xV'(2cx+b)dx}{\sqrt{1-x^2}}=1.
\]
Taking derivative according to $b$ we have
\[\int_{-1}^{1}\frac{V'(2cx+b)dx}{\sqrt{1-x^2}}=0.
\]
These equations will determine $b$ and $c.$ An $1-$cut potential $V$ must satisfying
\[1+c\int_{-1}^{1}{\frac{{V}'(2cy+b)\sqrt{1-{{y}^{2}}}dy}{y-x}}>0\] for $x\in (-1,1).$

\bigskip
\bigskip

{\footnotesize  
\par\noindent{\bf Acknowledgement.} 
The author would like to express his sincere thanks to Professor Academician Vilmos Totik for giving him the inversion problem of finite Hilbert transforms. 
Deepest appreciation is extended towards the NAFOSTED  (the National Foundation for Science and Techology Development in Vietnam) for the financial support.}

\end{document}